\documentclass[12pt]{article}
\usepackage{latexsym, amsfonts, amscd, amssymb, verbatim, amsmath,
enumerate}

\input{tcilatex}

\begin{document}

\title{On Radicals of Semirings and Related Problems}
\author{Y. Katsov \\
\textit{Department of Mathematics and Computer Science}\\
\textit{Hanover College, Hanover, IN 47243--0890, USA}\\
\textit{katsov@hanover.edu} \and T. G. Nam \\
\textit{Department of Mathematics}\\
\textit{Dong Thap Pedagogical University, Dong Thap, Vietnam}\\
\textit{trangiangnam05@yahoo.com}}
\date{}
\maketitle

\begin{abstract}
The aim of this paper is to develop an `external' Kurosh-Amitsur radical
theory of semirings and, using this approach, to obtain some fundamental
results regarding two Jacobson type of radicals --- the Jacobson-Bourne, $J$%
-, radical and a very natural its variation, $J_{s}$-radical --- of
hemirings, as well as the Brown-McCoy, $\mathcal{R}_{BM}$-, radical of
hemirings. Among the new central results of the paper, we single out the
following ones: Theorems unifying two, internal and external, approches to
the Kurosh-Amitzur radical theory of hemirings; A characterization of $J$%
-semisimple hemirings; A description of $J$-semisimple congruence-simple
hemirings; A characterization of finite additively-idempotent $J_{s}$%
-semisimple hemirings; Complete discriptions of $\mathcal{R}_{BM}$%
-semisimple commutative and lattice-ordered hemirings; Semiring versions of
the well-known classical ring results---Nakayama's and Hopkins Lemmas and
Jacobson-Chevalley Density Theorem; Establishing the fundamental
relationship between the radicals $J$, $J_{s}$, and $\mathcal{R}_{BM}$ of
hemirings $R$ and matrix hemirings $M_{n}(R)$; Establishing the
matric-extensibleness (see, \textit{e.g.}, \cite[Section 4.9]{gw:rtor}) of
the radical classes of the Jacobson, Brown-McCoy, and $J_{s}$-, radicals of
hemirings; Showing that the $J$-semisimplicity, $J_{s}$-semisimplicity, and $%
\mathcal{R}_{BM}$-semisimplicity of semirings are Morita invariant
properties.

\textbf{2010} \textbf{Mathematics Subject Classifications}: Primary 16Y60,
16D99, 16N99; Secondary 18A40, 18G05, 12K10

\end{abstract}

\section{Introduction}

As is well known, radical theory and radicals of algebraic structures,
originated in the beginning of the last century by Wedderburn and K\"{o}the
for rings (see, e.g., \cite{gw:rtor}), constitute important ``classical''
areas of the sustained interest in algebraic research which very often also
initiate research in new directions in other branches of mathematics. In
particular, studying of some analogs of the classical, Jacobson and
Brown-McCoy, ring radicals in a semiring setting commenced in 1950-60s in %
\cite{b:tjroas}, \cite{i:otjroas}, \cite{la:anotjroah}, and \cite{la:tbmroah}%
. As well, a Kurosh-Amitzur radical theory in a semiring setting has been
started in \cite{oj:rtfh}, \cite{on:anorfh}, \cite{ohl:wscoh}, \cite%
{olh:tscfh}, and \cite{ww:akartfps} (see, also, \cite{ww:ankrtfpsi} and \cite%
{ww:ankrtfpsii}), and then has been significantly advanced in \cite%
{hebwei:rtfs}, \cite{m:otrtfs}, \cite{hebwei:otibrtfsar}, and \cite%
{hebwei:scos}. Motivated by the Kurosh-Amitzur radical theory for rings
(see, e.g., \cite{gw:rtor}), the authors of the latter papers have developed
an `internal' Kurosh-Amitsur radical theory for semirings,\textit{\ i.e.},
the radical theory build ``axiomatically'' exclusively within the class $%
\mathbb{H}$ of all hemirings (semirings not necessarily possessing units)
and not involving representations of hemirings, that is, semimodules over
them. Although the main concepts of the both radical theories --- for rings
and hemirings --- are defined quite similarly, there are considerable
differences between these theories as well as all considerations and proofs
for semirings, not surprisingly, are significantly more complicated and
demand innovative ideas and techniques. As an algebraic objects, semirings
certainly are the most natural generalization of such algebraic systems as
rings and bounded distributive lattices, and investigating semirings and
their representations, one should undoubtedly use methods and techniques of\
both ring and lattice theory as well as diverse techniques and methods of
categorical and universal algebra. Thus, a wide variety of the algebraic
techniques used in studying semirings and their representations/semimodules
perhaps explains why research on a Kurosh-Amitsur radical theory for
semirings is still behind of that for rings.

The principal task here is to initiate an `external' Kurosh-Amitsur radical
theory for semirings --- a radical theory based, at this time, on
representations, semimodules, of semirings --- in the spirit of, and by
analogy with, the external radical theory for rings that can be found, for
example, in \cite[Section 3.14]{gw:rtor}; and then, based on it, to present
a series of fundamental results regarding the Jacobson and Brown-McCoy
radicals of hemirings, as well as to answer to some questions left open in
the earlier publications, mentioned above, on these topics.

The paper is organized as follows. In Section 2, for the reader's
convenience divided into two subsections, we included all subsequently
necessary notions and facts on semirings and semimodules, as well as on an
internal Kurosh-Amitzur radical theory of semirings.

In Section 3, after developing an external Kurosh-Amitsur radical theory of
semirings, we illustrate this approach by considering in more details two
Jacobson type of radicals --- the Jacobson-Bourne, $J$-, radical and a very
natural its variation, $J_{s}$-radical --- of hemirings, which coincide for
rings but are different in general, as well as the Brown-McCoy, $\mathcal{R}%
_{BM}$-, radical of hemirings. Among the new results of this section, which
are also among the central results of the paper, we single out the following
ones: Theorems 3.2, 3.3, and 3.4, actually unifying two, internal and
external, approaches to the Kurosh-Amitzur radical theory of hemirings;
Corollary 3.8, characterizing $J$-semisimple hemirings and answering a
question left open in \cite[Theorem 3.3, p. 12]{la:anotjroah}; a description
of $J$-semisimple congruence-simple hemirings (Theorem 3.10); a
characterization of finite additively-idempotent $J_{s}$-semisimple
hemirings (Theorem 3.11); complete descriptions of $\mathcal{R}_{BM}$%
-semisimple commutative and lattice-ordered hemirings (Corollaries 3.16 and
3.17).

In Section 4, reducing our semiring settings to corresponding\ ring ones, we
establish semiring versions of the well-known classical ring
results---Nakayama's and Hopkins Lemmas and Jacobson-Chevalley Density
Theorem (Theorem 4.3, Corollary 4.4, and Theorem 4.5, respectively), which
are among the main results of the paper.

Section 5 among other, in our view interesting and useful, observations,
contains the following main results of the paper: Theorem 5.8 and Corollary
5.11, establishing the fundamental relationship between the radicals $J$, $%
J_{s}$, and $\mathcal{R}_{BM}$ of hemirings $R$ and matrix hemirings $%
M_{n}(R)$, $n\geq 1$, and particularly extending \cite[Theorem 9]{b:tjroas}
from semirings to general hemirings; Theorem 5.14, establishing the
matric-extensibleness (see, \textit{e.g.}, \cite[Section 4.9]{gw:rtor}) of
the radical classes of the Jacobson, Brown-McCoy, and $J_{s}$-, radicals of
hemirings; and Theorem 5.17, showing that the $J$-semisimplicity, $J_{s}$%
-semisimplicity, and $\mathcal{R}_{BM}$-semisimplicity of semirings are
Morita invariant properties.

Finally, in the course of the paper, there have been stated several, in our
view interesting and promising, problems; also, all notions and facts of
categorical algebra, used here without any comments, can be found in \cite%
{macl:cwm}, and for notions and facts from semiring \ theory we refer to %
\cite{golan:sata}.

\section{Basic Concepts}

\subsection{Preliminaries on Semirings}

Recall \cite{golan:sata} that a \emph{hemiring\/} is an algebra $(R,+,\cdot
,0)$ such that the following conditions are satisfied:\medskip

(1) $(R,+,0)$ is a commutative monoid with identity element $0$;

(2) $(R,\cdot)$ is a semigroup;

(3) Multiplication distributes over addition from either side;

(4) $0r=0=r0$ for all $r\in R$.\medskip

And a hemiring $R$ is a \emph{semiring} if its multiplicative semigroup $%
(R,\cdot )$ actually is a monoid $(R,\cdot ,1)$ with identity element $1$.
Any hemiring $R$ can be naturally considered as an ideal of its \emph{Dorroh
extension} by $\mathbb{N}:=\{0,1,\ldots \}$ \cite[p.3]{golan:sata}, $%
R^{1}:=R\times \mathbb{N}$, which is a semiring with the identity element $%
(0,1)\in $ $R^{1}$ and operations of addition and multiplication defined by $%
(r,n)+(s,m):=(r+s,n+m)$ and $(r,n)(s,m):=(rs+mr+ns,nm)$ for all $%
(r,n),(s,m)\in R^{1}$.

A hemiring $R$ is \emph{additively cancellative} if $a+c=b+c$ implies $a=b$
for all $a,b,c\in R$. A nonempty subset $I$ of a hemiring $R$ is \emph{%
subtractive} if, for all $x,a\in R$, from $x+a,a\in I$ follows that $x\in I$%
, too. The \emph{subtractive closure} $\overline{I}$ of an ideal $I$ of a
hemiring $R$ is the smallest subtractive ideal of $R$ containing $I$, and it
is obviously defined as $\overline{I}:=$ $\{r\in R\,|\,r+x\in I$ for some $%
x\in I\}$; also, it is clear that an ideal $I$ a hemiring $R$ is subtractive
iff it coincides with its subtractive closure, \textit{i.e.}, $I=\overline{I}
$. Clearly, $0$ and $R$ are subtractive ideals for each hemiring $R$. By $%
\mathcal{I}(R)$ and $\mathcal{SI}(R)$ are denoted the sets of all ideals and
all subtractive ideals of a hemiring $R$, respectively; and a hemiring $R$
is called \emph{subtractive} if $\mathcal{I}(R)=\mathcal{SI}(R)$.

A hemiring $R$ is \emph{lattice-ordered} \cite[Section 21]{golan:sata} iff
it is also a lattice $(R,\vee ,\wedge )$, and for all $a,b\in R$, the
following conditions are satisfied: $a+b=a\vee b$ and, with respect to the
partial order naturally induced by the lattice operations, $ab\leq a\wedge b$%
.

As for rings, for any homomorphism $f:R\longrightarrow S$ between hemirings $%
R$ and $S$, there exists a subtractive ideal, the \emph{kernel}, $Ker(f):=$ $%
\{a\,|\,f(a)=0\}\subseteq R$ of $f$. A surjective hemiring homomorphism $%
f:R\longrightarrow S$ is a \emph{semiisomorphism} if $Ker(f)=0$. As usual,
the direct product $R=\prod_{i\in I}R_{i}$ of a family $(R_{i})_{i\in I}$ of
hemirings $R_{i}$ consists of the elements $r=(r_{i})_{i\in I}$ for $%
s_{i}\in R_{i}$ and is determined by the surjective homomorphisms $\pi
_{i}:R\longrightarrow R_{i}$ defined by $\pi _{i}(r)=r_{i}$; and a
subhemiring $S$ of $R$ is called a \emph{subdirect product} $S=\prod_{i\in
I}^{sub}R_{i}$ of $(R_{i})_{i\in I}$ if, for each $\pi _{i}$, the
restriction $\pi _{i}|_{S}:S\longrightarrow R_{i}$ is also surjective (see, 
\textit{e.g.}, \cite[p. 194]{hebwei:scos}).

Any ideal $I$ of a hemiring $R$ induces a congruence relation $\equiv _{I}$
on $R$, called the \emph{Bourne relation }\cite[p.78]{golan:sata}, given by $%
r\equiv _{I}r^{\prime }$ iff there exist elements $x,x^{\prime }\in I$ such
that $r+x=r^{\prime }+x^{\prime }$; and $R/I$ denotes the factor hemiring $%
R/\equiv _{I}$. Also, it is easy to see that the congruences $\equiv _{I}$
and $\ \equiv _{\overline{I}}$ on $R$ coincide for any ideal $I$ of $R$, and
hence, $R/I=R/\overline{I}$ hold for all ideals $I$ of $R$.

As usual, a \emph{left\/} $R$-\emph{semimodule} over a hemiring $R$ is a
commutative monoid $(M,+,0_{M})$ together with a scalar multiplication $%
(r,m)\mapsto rm$ from $R\times M$ to $M$ which satisfies the following
identities for all $r,r^{^{\prime }}\in R$ and $m,m^{^{\prime }}\in M$%
:\medskip

(1) $(rr^{^{\prime}})m=r(r^{^{\prime}}m)$;

(2) $r(m+m^{^{\prime }})=rm+rm^{^{\prime }}$;

(3) $(r+r^{^{\prime }})m=rm+r^{^{\prime }}m$;

(4) $r0_{M}=0_{M}=0m$.\medskip

\emph{Right semimodules\/} over a hemiring $R$ and homomorphisms between
semimodules are defined in the standard manner. If a hemiring $R$ is a
semiring, then all $R$-semimodules over $R$ are unitary ones. And, from now
on, let $\mathcal{M}_{R}$ and $_{R}\mathcal{M}$ denote the categories of
right and left semimodules, respectively, over a semiring $R$. As usual
(see, for example, \cite[Chapter 17]{golan:sata}), in the category $_{R}%
\mathcal{M}$ of left semimodules over a semiring $R$, a \textit{free} (left)
semimodule with a basis set $I$ is a direct sum (a coproduct) $\sum_{i\in
I}R_{i},R_{i}\cong $ $_{R}R$, $i\in I$, of $I$ copies of $_{R}R$; and a 
\textit{projective} left semimodule in $_{R}\mathcal{M}$ is just a retract
of a free semimodule. A semimodule $_{R}M$ is \textit{finitely generated}
iff it is a homomorphic image of a free semimodule with a finite basis set.
A left semimodule $M$ over a hemiring $R$ is \emph{cancellative} if $x+z=y+z$
implies $x=y$ for all $x,y,z\in M$. A nonempty subset $N$ of an $R$%
-semimodule $M$ is \emph{subtractive} if, for all $x,y\in M$, from $x+y,x\in
N$ follows that $y\in N$, too.

The usual concepts of the \textit{Descending}\textit{\ Chain Condition} and 
\textit{artinian} modules of theory of modules over rings, as well as
results involving them, are easily extended in an obvious\ fashion (see, 
\textit{e.g.}, \cite{knt:ossss}) to a context of semimodules over semirings.
As for the modules over rings, $(0:M)$ or $(0:M)_{R}$ denotes the \emph{%
annihilator }of a left $R$-semimodule $_{R}M$, \textit{i.e.}, $%
(0:M)_{R}:=\{r\in R\,|\,rM=0\}$; and $_{R}M$ is \emph{faithful} iff $%
(0:M)_{R}=0$. Finally and similarly to the case of module over rings, the
following observations can be easily verified and will prove to be useful in
sequel:\smallskip

\noindent \textbf{Proposition 2.1. }(\textit{cf.} \cite[Proposition 3.14.1]%
{gw:rtor}) \textit{Let }$R$\textit{\ be a hemiring, and }$I$\textit{\ }$\in $
$\mathcal{I}(R)$\textit{. }

(i) \textit{If }$M$ \textit{is a left }$R/I$\textit{-semimodule, then under
the scalar multiplication} 
\begin{equation*}
rm=\overline{r}m\text{,}
\end{equation*}%
$M$ \textit{becomes a left} $R$\textit{-semimodule with} $I\subseteq
(0:M)_{R}$.

(ii) \textit{If }$M$ \textit{is a left }$R$\textit{-semimodule with} $%
I\subseteq (0:M)_{R}$, \textit{then} $M$ \textit{is a left }$R/I$\textit{%
-semimodule under the scalar multiplication} 
\begin{equation*}
\overline{r}m=rm.
\end{equation*}

(iii) \textit{Every subsemimodule of the left} $R/I$\textit{-semimodule} $M$ 
\textit{is a subsemimodule of the left} $R$\textit{-semimodule} $M$\textit{,
and the converse is also true when} $I\subseteq (0:M)_{R}$.

(iv) $(0: M)_{R/I} = (0: M)_R/I$\textit{.\medskip}

\noindent \textbf{Proof.} (i), (ii) and (iii) are clear.

(iv) Obliviously, $(0:M)_{R/I}\supseteq (0:M)_{R}/I$. Conversely, for any $%
\overline{r}\in (0:M)_{R/I}$, $\overline{r}m=0$ for all $m\in M$, and,
hence, $rm=0$ for all $m\in M.$ Thus $r\in (0:M)_{R}$. $\ \ \ \ _{\square
}\medskip $ \ \ \ \textit{\ \ \ \ \ \ }

\subsection{Preliminaries on Radical Theory of Semirings}

In this subsection, we briefly sketch the basic concepts of an `internal'
Kurosh-Amitsur radical theory of semirings---in other words, the
Kurosh-Amitsur radical theory build exclusively within the class $\mathbb{H}$
of all hemirings without using representations of hemirings---originated by
D. M. Olson and his coauthors in a series of papers \cite{oj:rtfh}--\cite%
{olh:tscfh}, and then considerably advanced in \cite{hebwei:rtfs}, \cite%
{m:otrtfs}, \cite{hebwei:otibrtfsar}, and \cite{hebwei:scos}. As was shown
in the latter papers (see, for instance, \cite{m:otrtfs}), similarly to the
radical theory of rings\ (see, \textit{e.g.}, \cite{gw:rtor}), there are
three equivalent approaches to the Kurosh-Amitsur radical theory of
semirings---by means of radical classes, radical operators, and semisimple
classes, independently defined in a fixed universal class $\mathbb{U}$ $%
\subseteq \mathbb{H}$ of hemirings. Herewith, a nonempty subclass $\mathbb{U}
$ of $\mathbb{H}$ is called \emph{universal} if $\mathbb{U}$ is a \emph{%
hereditary} ($R\in \mathbb{U}$ implies $\mathcal{I}(R)\subseteq \mathbb{U}$)
and \emph{homomorphically closed} ($R\in \mathbb{U}$ implies $\varphi (R)\in 
\mathbb{U}$ for every hemiring homomorphism $\varphi $) subclass in $\mathbb{%
H}$. Furthermore, by \cite[Definition 3.1]{m:otrtfs}, a nonempty subclass $%
\mathbb{R}$ of a fixed universal class $\mathbb{U}\subseteq $ $\mathbb{H}$
is called a \emph{radical class} of $\mathbb{U}$ if $\mathbb{R}$ satisfies
the following two conditions: (i) $\mathbb{R}$ is homomorphically closed;
(ii) For every hemiring $R\in \mathbb{U}\setminus \mathbb{R}$, there exists
a subtractive ideal $I\in \mathcal{SI}(R)\setminus \{R\}$ such that $%
\mathcal{I}(R/I)\cap \mathbb{R}=\{(0)\}$, where $(0)$ is the trivial, zero,
hemiring.

Analogously to the case of rings (see, \textit{e.g.}, \cite[Chapter II]%
{gw:rtor}), there have been established the following characterizations of
radical classes for hemirings:\smallskip

\noindent \textbf{Theorem 2.2.} (\cite[Theorems 3.2 and 4.7]{m:otrtfs}) 
\textit{For a subclass} $\mathbb{R}$ \textit{of} \textit{a universal class} $%
\mathbb{U}$\textit{,} \textit{the following conditions are equivalent:}

(1) $\mathbb{R}$ \textit{is a radical class of} $\mathbb{U}$\textit{.}

(2) $\mathbb{R}$ \textit{satisfies the following two properties:}

\quad (i) \textit{If} $R\in \mathbb{R}$\textit{,} \textit{then for every
nonzero surjective hemiring homomorphism }$R\longrightarrow S$ \textit{there
exists a nonzero ideal} $I$ \textit{of} $S$ \textit{such that} $I\in \mathbb{%
R}$\textit{;}

\quad (ii) \textit{If} $R\in \mathbb{U}$ \textit{and for every nonzero
surjective hemiring homomorphism} $R\longrightarrow S$\textit{\ there exists
a nonzero ideal} $I$ \textit{of} $S$ such that $I\in \mathbb{R}$\textit{,} 
\textit{then} $R\in \mathbb{R}$.

(3) \textit{The following three properties are true for }$\mathbb{R}$\textit{%
:}

\quad (i) $\mathbb{R}$ \textit{is homomorphically closed;}

\quad (ii) $\mathbb{R}$ \textit{is extensionally closed in} $\mathbb{U}$ 
\textit{--- for} \textit{all} $R\in \mathbb{U}$ \textit{and }$I\in \mathcal{I%
}(R)$\textit{, if }$I\in \mathbb{R}$ \textit{and} $R/I\in \mathbb{R}$\textit{%
,} \textit{then} $R\in \mathbb{R}$\textit{;}

\quad (iii) $\mathbb{R}$ \textit{has the inductive property --- for any }$%
R\in \mathbb{U}\ $\textit{and ascending chain of ideals} $I_{1}\subseteq
...\subseteq I_{\lambda }\subseteq ...$\textit{with all} $I_{\lambda }\in 
\mathcal{I}(R)\cap \mathbb{R}$\textit{,} \textit{the ideal} $\cup I_{\lambda
}\in \mathbb{R}$\textit{,} \textit{too.}$\medskip $

As was shown in \cite[Theorem 5.3]{m:otrtfs} (and for a more restricted
class of universal classes in \cite[p. 309 and Theorem 17]{on:anorfh}, too),
the general method of a construction of the \textit{upper radical} $\mathcal{%
U}\mathbb{S}$ for a regular class of rings $\mathbb{S}$ (see, e.g., %
\cite[Theorem 2.2..3]{gw:rtor}) can be extended to a hemiring setting as
well. Recall (see \cite[p. 542 and Theorem 5.3]{m:otrtfs} or \cite[p. 182]%
{hebwei:otibrtfsar}) that a subclass $\mathbb{S}$ of a universal class $%
\mathbb{U}$ $\subseteq \mathbb{H}$ is a \emph{regular class} if, for any
nonzero ideal $I$ of a hemiring $R\in \mathbb{S}$, there exists a nonzero
surjective hemiring homomorphism $I\longrightarrow S$ with $S\in \mathbb{S}$%
. Since Theorem 5.3 of \cite{m:otrtfs} was proved under the assumption that
a regular subclass $\mathbb{S}$ of a universal class $\mathbb{U}$ is an
isomorphically closed class, we have found to be reasonable to propose here
an alternative, new proof of it eliminating this assumption. \smallskip

\noindent \textbf{Theorem 2.3. }(\textit{cf.} \cite[Theorem 5.3]{m:otrtfs}) 
\textit{If} $\mathbb{S}$ \textit{is a regular subclass of a universal class} 
$\mathbb{U}$ $\subseteq \mathbb{H}$, \textit{then the class}%
\begin{equation*}
\mathcal{U}\mathbb{S}=\{R\in \mathbb{U}\,|\,\text{\ }R\text{\textit{has no
nonzero homomorphic image in} }\mathbb{S}\}
\end{equation*}%
\newline

\noindent \textit{is a radical class of} $\mathbb{U}$ \textit{and} $\mathbb{S%
}\cap \mathbb{U}\mathbb{S}=\{(0)\}$\textit{; moreover, }$\mathcal{U}\mathbb{S%
}$ \textit{is the largest radical class in }$\mathbb{U}\,$\textit{having
zero intersection with} $\mathbb{S}$\textit{.}$\medskip $

\noindent \textbf{Proof.} It is obvious, as in the proof of \cite[Theorem
5.3]{m:otrtfs}, that $\mathcal{U}\mathbb{S}$ is homomorphically closed.

Let us show that $\mathcal{U}\mathbb{S}$ is extensionally closed. Suppose
that for a hemiring $R$ and and ideal $I\in \mathcal{I}(R)$ we have that $%
I,R/I\in \mathcal{U}\mathbb{S}$, but $R\notin \mathcal{U}\mathbb{S}$. Then,
there exists a nonzero surjective homomorphism $f:R\longrightarrow S$ with $%
S\in \mathbb{S}$, and consider two possible cases:

a) $I\nsubseteq Ker(f)$. Then, for $f(I)$ is a nonzero ideal of $S\in 
\mathbb{S}$ and $\mathbb{S}$ is a regular class, there is a nonzero
surjective homomorphism $g:f(I)\longrightarrow H$ with $H\in \mathbb{S}$;
hence, $gf|_{I}:I\longrightarrow H$, where $f|_{I}$ \ is the restriction of $%
f$ on $I$, is a nonzero surjective homomorphism. Therefore, one gets a
contradiction $I\notin \mathcal{U}\mathbb{S}$.

b) $I\subseteq Ker(f)$. Then, for the mapping $h:R/I\longrightarrow S$,
defined by $h(\overline{r})=f(r)$, is a nonzero surjective homomorphism, one
has a contradiction $R/I\notin \mathcal{U}\mathbb{S}$.

Let us show that $\mathcal{U}\mathbb{S}$ has the inductive property. Suppose
that for $R\in \mathbb{U}$, a chain $I_{1}\subseteq ...\subseteq I_{\lambda
}\subseteq ...$ with all $I_{\lambda }\in \mathcal{I}(R)\cap \mathcal{U}%
\mathbb{S}$, and $I=\cup I_{\lambda }$, we have $I\notin \mathcal{U}\mathbb{S%
}$. Then, there exists a nonzero surjective homomorphism $f:I\longrightarrow
J$ with $J\in \mathbb{S}$. Hence, there exists an element $a_{\lambda }\in
I_{\lambda }\subseteq I$ such that $f(a_{\lambda })\neq 0$, and $%
f|_{I_{\lambda }}(I_{\lambda })$ is a nonzero ideal of $J\in \mathbb{S}$.
Since $\mathbb{S}$ is a regular class, there exist nonzero surjective
homomorphisms $g:f|_{I_{\lambda }}(I_{\lambda })\longrightarrow K$ $\in 
\mathbb{S}$ and $gf|_{I_{\lambda }}:I_{\lambda }\longrightarrow K$ $\in 
\mathbb{S}$, and, therefore, we have a contradiction $I_{\lambda }\notin 
\mathcal{U}\mathbb{S}$.

By Theorem 2.2 (3), we have established that $\mathcal{U}\mathbb{S}$ is a
radical class of $\mathbb{U}$.

The rest is almost obvious and can be shown as in \cite[Theorem 5.3]%
{m:otrtfs}. $\ \ \ \ _{\square }\medskip $

Following \cite[Definition 2.4]{hebwei:otibrtfsar} (or \cite[Definition 4.1]%
{m:otrtfs}), a mapping $\varrho :\mathbb{U}\longrightarrow \mathbb{U}$ is
called a \emph{radical operator} in $\mathbb{U}$ if each hemiring $R\in 
\mathbb{U}$ has an image $\varrho (R)\in \mathcal{SI}(R)\subseteq \mathbb{U}$
such that the following conditions are satisfied for all $R,S\in \mathbb{U}$:

(i) $\varphi (\varrho (R))\subseteq \varrho (\varphi (R))$ for each
homomorphism $\varphi :R\longrightarrow S$\textit{;}

(ii) $\varrho (R/\varrho (R))=(0)$\textit{;}

(iii) \textit{For every nonzero ideal }$I$\textit{\ of }$R$\textit{, }$%
\varrho (I)=I$\textit{\ implies that }$I\subseteq \varrho (R)$\textit{;}

(iv) $\varrho (\varrho (R))=\varrho (R)$.\smallskip

Terming an ideal $I\in \mathcal{I}(R)\cap \mathbb{R}$ of a hemiring $R$ for
a radical class $\mathbb{R}$ an $\mathbb{R}$\textit{-ideal} of $R$, it was
shown \cite[Theorem 3.7]{m:otrtfs} that all $\mathbb{R}$-ideals of $R\in 
\mathbb{U}$ are contained in the greatest ideal of this kind, called the $%
\mathbb{R}$\textit{-radical} of $R$. More precisely, we have the following
fact:\medskip

\noindent \textbf{Theorem 2.4.} (\cite[Theorem 2.3]{hebwei:otibrtfsar}) 
\textit{For a radical class} $\mathbb{R}$ of $\mathbb{U}$ \textit{and} $R\in 
\mathbb{U}$\textit{,} \textit{the union} $\varrho _{\mathbb{R}}(R)$ \textit{%
of all} $\mathbb{R}$\textit{-ideals of} $R$ \textit{is again an} $\mathbb{R}$%
\textit{-ideal of} $R$\textit{,} \textit{i.e.,} $\varrho _{\mathbb{R}%
}(R)=\bigcup \{I\in \mathcal{I}(R)\cap \mathbb{R}\}\in \mathcal{I}(R)\cap 
\mathbb{R}.$ \textit{Moreover,} $\varrho _{\mathbb{R}}(R)$ \textit{is a
subtractive ideal of} $R$\textit{,} \textit{and the mapping} $\varrho _{%
\mathbb{R}}:\mathbb{U}\longrightarrow \mathbb{U}$ \textit{determined in this
way is a radical operator in} $\mathbb{U}$\textit{.}$\medskip $

A subclass $\mathbb{S}$ of a universal class $\mathbb{U}$ is called a \emph{%
semisimple class} of $\mathbb{U}$ \cite[Definition 2.5]{hebwei:otibrtfsar}
iff $\mathbb{S}$ satisfies the following two conditions:

(i) If $R\in \mathbb{S}$, then for every nonzero ideal $I$ of $R$ there
exists a nonzero surjective hemiring homomorphism $I\longrightarrow S$ onto $%
S\in \mathbb{S}$;

(ii) If $R\in \mathbb{U}$ and for every nonzero ideal $I$ of $R$ there
exists a nonzero surjective hemiring homomorphism $I\longrightarrow S$ onto $%
S\in \mathbb{S}$, then $R\in \mathbb{S}$.\smallskip

As has been shown in \cite[Theorems 4.5, 5.4, and 5.5]{m:otrtfs} (see also %
\cite[Theorem 2.6]{hebwei:otibrtfsar}), each of the three concepts
above---radical classes, radical operators, and semisimple classes---can
serve as a starting point for an internal Kurosh-Amitsur radical theory of
semirings/hemirings; and, of course, a hemiring is semisimple iff its
radical is zero.

\section{Semimodules and Radical Classes}

We start this section with a developing of an `external' Kurosh-Amitsur
radical theory for semirings/hemirings---a Kurosh-Amitsur radical theory
based, at this time, on representations, semimodules, of hemirings---in the
spirit of, and by analogy with, the external radical theory for rings (see,
for example, \cite[Section 3.14]{gw:rtor}). Hear again, as was mentioned
earlier, despite of the natural fact that many corresponding concepts and
considerations for radical theories for rings and semirings/hemirings are
very similar, in the semiring setting, they are, not surprisingly,
significantly more complicated and often---since there is no ``pleasure'' of
the abelian categories, and one must ``live in the world without
subtraction''---involve original ideas and techniques.

Now, for each hemiring $R$, let $\Sigma _{R}$ be a ``chosen'' class of left $%
R$-semimodule $M$ with $RM\neq 0$, and $\Sigma $ be the union of all those $%
\Sigma _{R}$. Let 
\begin{equation*}
ker(\Sigma _{R}):=\cap \{(0:M)_{R}\,|\,M\in \Sigma _{R}\}.
\end{equation*}

Then, let us single out the following four conditions on the class of
semimodules $\Sigma $ which it might satisfy:

(SM$1$) If $M\in \Sigma _{R/I}$, then $M\in \Sigma _{R}$;

(SM$2$) If $M\in \Sigma _{R}$ and $I$ is an ideal of $R$ such that $%
I\subseteq (0:M)_{R}$, then $M\in \Sigma _{R/I}$;

(SM$3$) If $ker(\Sigma _{R})=0$, then $\Sigma _{I}\neq \emptyset $ for all
nonzero ideals $I$ of $R$;

(SM$4$) If $\Sigma _{I}\neq \emptyset $ for each nonzero ideal $I$ of $R$,
then $ker(\Sigma _{R})=0$.\medskip

Denoting by $\mathcal{F}(\Sigma )$ the class of hemirings $R$ for which the
semimodule class $\Sigma _{R}$ contains a faithful semimodule, using
Proposition 2.1 and repeating verbatim the proof of \cite[Proposition 3.14.2]%
{gw:rtor}, one obtains its hemiring analog:\medskip

\noindent \textbf{Proposition 3.1. }(\textit{cf. }\cite[Proposition 3.14.2]%
{gw:rtor}) $\mathcal{F}(\Sigma )$ \textit{is a regular class.\smallskip }

The following three observations are hemiring versions of the corresponding
ring results.\smallskip

\noindent \textbf{Theorem 3.2.} (\textit{cf. }\cite[Theorem 3.14.3]{gw:rtor}%
) \textit{Let} $\Sigma $ \textit{be a class of semimodules satisfying
conditions} (SM$1$)$-$(SM$3$), \textit{and }$\Re (\Sigma ):=\{R\in \mathbb{H}%
\,|\,\Sigma _{R}=\emptyset \}$. \textit{Then:}

\textit{(i)} $\Re (\Sigma )$ \textit{is a radical class, and} $\Re (\Sigma )=%
\mathcal{U}\mathcal{F}(\Sigma )$\textit{;}

\textit{(ii)} $\varrho _{\Re (\Sigma )}(R)\subseteq ker(\Sigma _{R})$ 
\textit{for every hemiring} $R$\textit{.\medskip }

\noindent \textbf{Proof.} (i) This can be justified by applying\ Proposition
2.1 and Theorem 2.3 and repeating verbatim the proof of \cite[Theorem 3.14.3
(i)]{gw:rtor}.

(ii) For every hemiring $R$ and semimodule $M\in \Sigma _{R}$, the
annihilator\emph{\ }$(0:M)$ is obviously a subtractive ideal of $R$, and
therefore, by \cite[Theorem 2.2]{m:otrtfs}, there is the semiisomorphism $%
f:\varrho _{\Re (\Sigma )}(R)/(\varrho _{\Re (\Sigma )}(R)\cap
(0:M))\longrightarrow (\varrho _{\Re (\Sigma )}(R)+(0:M))/(0:M)$ given by $%
\overline{r}\longmapsto \overline{r}$.

It is clear that $(\varrho _{\Re (\Sigma )}(R)+(0:M))/(0:M)$ is an ideal of $%
R/(0:M)\in \mathcal{F}(\Sigma )$. Assuming that $(\varrho _{\Re (\Sigma
)}(R)+(0:M))/(0:M)\neq 0$ and using Proposition 2.1 and Theorem 2.3, one has
nonzero surjective homomorphisms $g:(\varrho _{\Re (\Sigma
)}(R)+(0:M))/(0:M)\longrightarrow S$ $\in \mathcal{F}(\Sigma )$ and $gf\pi
:\varrho _{\Re (\Sigma )}(R)\longrightarrow S$ $\in \mathcal{F}(\Sigma )$,
where $\pi :\varrho _{\Re (\Sigma )}(R)\longrightarrow \varrho _{\Re (\Sigma
)}(R)/(\varrho _{\Re (\Sigma )}(R)\cap (0:M))$ is the natural surjection.
The latter implies that, in contradiction with (i), $\varrho _{\Re (\Sigma
)}(R)\notin \mathcal{U}\mathcal{F}(\Sigma )$. Hence, $(\varrho _{\Re (\Sigma
)}(R)+(0:M))/(0:M)=0$ and $\varrho _{\Re (\Sigma )}(R)+(0:M)=(0:M)$; and,
therefore, $\varrho _{\Re (\Sigma )}(R)\subseteq (0:M)$ and $\varrho _{\Re
(\Sigma )}(R)\subseteq ker(\Sigma _{R})$. \textit{\ \ \ \ \ \ }$_{\square
}\medskip $

\noindent \textbf{Theorem 3.3.} (\textit{cf. }\cite[Theorem 3.14.4]{gw:rtor}%
) \textit{Let} $\Sigma $ \textit{be a class of semimodules satisfying
conditions} (SM$1$)--(SM$4$). \textit{Then:}

\textit{(i)} $R$ \textit{is} $\Re (\Sigma )$\textit{-semisimple iff} $R$ 
\textit{is semiisomorphic to a subdirect product of hemirings in} $\mathcal{F%
}(\Sigma )$\textit{;}

\textit{(ii)} $\varrho _{\Re (\Sigma )}(R)=ker(\Sigma _{R})$ \textit{for
every hemiring} $R$.\smallskip

\noindent \textbf{Proof.} As was shown in Theorem 3.2, $\Re (\Sigma )$ is a
radical class, and, hence, by \cite[Theorem 2.6(c)]{hebwei:otibrtfsar}, the
class $\mathbb{S}:=\{R\in \mathbb{H}\,|\,\varrho _{\Re (\Sigma )}(R)=(0)\}$
is a semisimple class of $\mathbb{H}$.

(i). $\Longrightarrow $. Let $R$ be $\Re (\Sigma )$-semisimple. Then, $\
\varrho _{\Re (\Sigma )}(R)=0$, and $R$ has no nonzero ideals in $\Re
(\Sigma )$, and hence, $\Sigma _{I}\neq \emptyset $ for all nonzero ideals $%
I $ of $R$. From the latter and condition (SM$4$), we have $ker(\Sigma
_{R})=\cap _{M\in \Sigma _{R}}(0:M)=0$. By Proposition 2.1, $R/(0:M)\in 
\mathcal{F}(\Sigma )$ for all $M\in \Sigma _{R}$; and, by \cite[Theorem 6.3
(b)]{hebwei:scos}, there exists the semiisomorphism $f:R\longrightarrow
\prod_{M\in \Sigma _{R}}R/(0:M)$ given by $r\longmapsto ([r]_{\equiv
_{(0:M)}})_{M\in \Sigma _{R}}$.

$\Longleftarrow $. Let $R$ be semiisomorphic to a subdirect product $%
T=\prod_{i\in I}^{sub}R_{i}$ of hemirings $R_{i}\in \mathcal{F}(\Sigma ).$
Then, $ker(\Sigma _{R_{i}})=0$, and hence, by Theorem 3.2 (ii), $R_{i}$ is $%
\Re (\Sigma )$-semisimple, and $R_{i}\in \mathbb{S}$ for all $i\in I$. By %
\cite[Theorem 4.3 (a)]{hebwei:scos} $\mathbb{S}$ is subdirectly closed in $%
\mathbb{H}$ and, therefore, $T=\prod_{i\in I}^{sub}R_{i}\in \mathbb{S}$,
too. By \cite[Theorem 3.7(b)]{hebwei:scos}, $\mathbb{S}$ is also inverse
semiisomorphically closed in $\mathbb{H}$ (\textit{i.e.}, if there exists a
semiisomorphism $\varphi :R\rightarrow T$ $\in \mathbb{S}$, then $R\in 
\mathbb{S}$) and, hence, $R\in \mathbb{S}$, \textit{i.e.}, $R$ is $\Re
(\Sigma )$-semisimple.

(ii). Using Theorem 3.2 (ii), we only need to show that $\varrho _{\Re
(\Sigma )}(R)\supseteq ker(\Sigma _{R})$. Indeed, since $R/\varrho _{\Re
(\Sigma )}(R)$ is clearly $\Re (\Sigma )$-semisimple, as in the proof of the
first part of (i), we get the semiisomorphism $f:R/\varrho _{\Re (\Sigma
)}(R)\longrightarrow \prod_{M\in \Sigma _{R/\varrho _{\Re (\Sigma
)}(R)}}^{sub}R/\varrho _{\Re (\Sigma )}(R)/(0:M)_{R/\varrho _{\Re (\Sigma
)}(R)}$ with

\noindent $R/\varrho _{\Re (\Sigma )}(R)/(0:M)_{R/\varrho _{\Re (\Sigma
)}(R)}\in \mathcal{F}(\Sigma )$; and, hence, $R/\varrho _{\Re (\Sigma )}(R)$
is semiisomorphic to a subdirect product $f(R/\varrho _{\Re (\Sigma )}(R))$
of hemirings $R/\varrho _{\Re (\Sigma )}(R)/(0:M)_{R/\varrho _{\Re (\Sigma
)}(R)}$. Applying Proposition 2.1(iv) and \cite[Theorem 2.3]{hebwei:scos},
one gets

\noindent $(0:M)_{R/\varrho _{\Re (\Sigma )}(R)}\cong (0:M)_{R}/\varrho
_{\Re (\Sigma )}(R)$ and $R/\varrho _{\Re (\Sigma )}(R)/(0:M)_{R/\varrho
_{\Re (\Sigma )}(R)}\cong R/(0:M)_{R}$. Hence, $R/(0:M)_{R}\in \mathcal{F}%
(\Sigma )$, and therefore there is a faithful left $R/(0:M)_{R}$-semimodule $%
N$ in $\Sigma _{R/(0:M)_{R}}$; and for condition (SM$1$), $N\in \Sigma _{R}$%
. Now, applying again Proposition 2.1, we get $%
(0:N)_{R}/(0:M)_{R}=(0:N)_{R/(0:M)_{R}}=0$, and hence, $ker(\Sigma
_{R})\subseteq (0:N)_{R}=(0:M)_{R}$ and $ker(\Sigma _{R})\subseteq \cap
(0:M)_{R}\subseteq \varrho _{\Re (\Sigma )}(R)$.\textit{\ \ \ \ \ \ }$%
_{\square }\medskip $

The next observation, actually inverse to Theorems 3.2 and 3.3, shows that
every radical class can be obtained from a suitable class of semimodules $%
\Sigma $ satisfying conditions (SM$1$)$-$(SM$4$).\smallskip

\noindent \textbf{Theorem 3.4.} (\textit{cf. }\cite[Theorem 3.14.5]{gw:rtor}%
) \textit{Let} $\mathbb{R}$ $\subseteq \mathbb{H}$ \textit{be a radical
class, }$\Sigma _{R}=\{_{R}M\in |_{R}\mathcal{M}|\,|\,RM\neq 0$ $\mathit{\&}$%
\textit{\ }$\varrho _{\mathbb{R}}(R/(0:M))=(0)\}$ \textit{for every} $R$ $%
\in $ $\mathbb{H}$\textit{,} \textit{and} $\Sigma =\cup \Sigma _{R}$\textit{%
. Then:}

\textit{(i)} $\Sigma $ \textit{satisfies conditions} (SM$1$)--(SM$4$)\textit{%
;}

\textit{(ii)} $\mathbb{R}=\Re (\Sigma )$\textit{.\medskip }

\noindent \textbf{Proof.} (i) Using Proposition 2.1 and repeating verbatim
the proofs for the corresponding ring conditions (M$1$)$\ $and (M$2$) in %
\cite[Theorem 3.14.5]{gw:rtor}, we get that the class $\Sigma $ satisfies (SM%
$1$)$\ $and (SM$2$).

To establish (SM$3$) for the class $\Sigma $, suppose that $ker(\Sigma
_{R})=0$ and $I$ is a nonzero ideal of $R$, and therefore, $IM\neq 0$ for
some semimodule $M\in \Sigma _{R}$. By \cite[Theorem 2.6(c)]%
{hebwei:otibrtfsar}, the class of semimodules $\mathbb{S}=\{R\in \mathbb{H}%
\,|\,\varrho _{\mathbb{R}}(R)=(0)\}$ is a semisimple class of $\mathbb{H}$.
Then, by \cite[Theorem 2.2]{m:otrtfs} there exists a semiisomorphism $%
f:I/(0:M)_{I}=I/((0:M)_{R}\cap I)\longrightarrow (I+(0:M)_{R})/(0:M)_{R}$,
where $(I+(0:M)_{R})/(0:M)_{R}$ is an ideal of $R/(0:M)_{R}\in \mathbb{S}$.
Whence, by \cite[Theorem 3.7(a)]{hebwei:scos} $(I+(0:M)_{R})/(0:M)_{R}\in 
\mathbb{S}$, and since by \cite[Theorem 3.7(b)]{hebwei:scos} $\mathbb{S}$ is
inverse semiisomorphically closed in $\mathbb{H}$, we have $I/(0:M)_{I}\in 
\mathbb{S}$ and, hence, $\varrho _{\mathbb{R}}(I/(0:M)_{I})=(0)$. From the
latter, one has that $M\in \Sigma _{I}$ and, therefore, $\Sigma $ satisfies
(SM$3$), too.

The fact that $\Sigma $ satisfies (SM$4$) can be establish by verbatim
repeating the proof of the condition (M$4$) for $\Sigma $ in \cite[Theorem
3.14.5]{gw:rtor} and using Theorem 2.2.

(ii) Suppose $R\in \mathbb{R}$ and $M$ is a left $R$-semimodule having $%
\varrho _{\mathbb{R}}(R/(0:M))=(0)$. Then, by Theorem 2.2 $R/(0:M)\in 
\mathbb{R}$ and, hence, $R/(0:M)=\varrho _{\mathbb{R}}(R/(0:M))=(0)$.
Therefore, as $(0:M)$ is a subtractive ideal of $R$, we have $R=(0:M)$ and,
hence, $RM=0$. From the latter we conclude that $\Sigma _{R}=\emptyset $,
and, therefore, $R\in \Re (\Sigma )$, \textit{i.e.}, $\mathbb{R}\subseteq
\Re (\Sigma )$.

Now suppose $R\notin \mathbb{R}$. From Theorem 2.2 it follows that there
exists a nonzero surjective homomorphism $f:R\longrightarrow S$ with $%
\varrho _{\mathbb{R}}(S)=(0)$, and, therefore, $S$ is semiisomorphic $R/I$,
where $I=Ker(f)$, and applying \cite[Theorem 3.7(b)]{hebwei:scos}, one gets $%
\varrho _{\mathbb{R}}(R/I)=(0)$. Then, consider the Dorroh extension $%
(R/I)^{1}$of $R/I$: It is obvious that $(R/I)^{1}$ is a left $R/I$%
-semimodule with $(0:(R/I)^{1})_{R/I}=(0)$; hence, $%
(R/I)/(0:(R/I)^{1})_{R/I}=R/I$ is $\mathbb{R}$-semisimple. The latter
implies that $(R/I)^{1}\in \Sigma _{R/I}$; moreover, since $\Sigma $
satisfies (SM$1$), $(R/I)^{1}\in \Sigma _{R}$, and, hence, $R\notin \Re
(\Sigma )$. Thus, $\Re (\Sigma )\subseteq \mathbb{R}$. \textit{\ \ \ \ \ \ }$%
_{\square }\medskip $

As Theorem 3.4 shows, any radical class might be, in principle,
``externally'' described by means of the class of suitable
presentations/semimodules of hemirings. However, in this theorem the class
of the corresponding semimodules $\Sigma $, constructed straightforwardly
from the radical itself, is the greatest of all possible classes describing
a given radical class. Therefore, it is more important and interesting to be
able to find a smaller than $\Sigma $ class of semimodules describing the
same radical class and/or having a nice practical characterization. In some
important cases this can be successfully done and to illustrate that, in the
following two subsections of this section, we consider two analogs of the
classical ring radicals---the Jacobson and the Brown-McCoy radicals---for
semirings.

\subsection{The Jacobson Type Radicals of Semirings}

A nonzero cancellative left semimodule $M$ over a hemiring $R$ is \emph{%
irreducible }\cite[Definition 5]{i:otjroas} if, for an arbitrarily fixed
pair of elements $m_{1},m_{2}$ of $M$ with $m_{1}\neq m_{2}$ and any $m\in M$%
, there exist $r_{1},r_{2}\in R$ such that 
\begin{equation*}
m+r_{1}m_{1}+r_{2}m_{2}=r_{1}m_{2}+r_{2}m_{1}\text{.}
\end{equation*}

It is easy to see that $RM\neq 0$ for any irreducible left $R$-semimodule $M$%
. And a nonzero cancellative left semimodule $M$ over a hemiring $R$ is 
\emph{semi-irreducible} \cite[p. 412]{i:otjroas} if $RM\neq 0$ and there are
no subtractive subsemimodules of $M$ except $M$ itself and the zero one.
Finally, we call a left $R$-semimodule $M$ \emph{simple} if $RM\neq 0$ and
there are only trivial subsemimodules of, as well as congruences on, $M$%
.\medskip

\noindent \textbf{Proposition 3.5. }(\textit{cf.} \cite[Theorem 3.14.6]%
{gw:rtor}) \textit{For every hemiring} $R$ $\in $ $\mathbb{H}$\textit{, let} 
$\Sigma _{R}$ \textit{denote the corresponding class of semimodules in each
of the cases:}

\textit{(i)} $\Sigma _{R}:=\{_{R}M\in |_{R}\mathcal{M}|\,|\,M$ \textit{is
irreducible}$\}$\textit{,}

\textit{(ii)} $\Sigma _{R}:=\{_{R}M\in |_{R}\mathcal{M}|\,|\,M$ \textit{is
semi-irreducible}$\}$\textit{,}

\textit{(iii)} $\Sigma _{R}:=\{_{R}M\in |_{R}\mathcal{M}|\,|\,M$ \textit{is
simple}$\}$\textit{.}

\noindent \textit{Then,} \textit{for each of these cases,} $\Sigma :=\cup
\Sigma _{R}\ $\textit{satisfies conditions} (SM$1$)--(SM$4$)\textit{%
.\medskip }

\noindent \textbf{Proof.} First, for each of\ the cases (i)--(iii), one may
easily verify that $\Sigma $ satisfies (SM$1$) and (SM$2$).

(i). Suppose $ker(\Sigma _{R})=0$ for a hemiring $R$, and $I$ is a nonzero
ideal of $R$. Then, $IM\neq 0$ for some irreducible left $R$-semimodule $M$
and by \cite[Lemma 4(1)]{i:otjroas} $M$ is an irreducible left $I$%
-semimodule; so, $\Sigma _{I}\neq \emptyset $ and $\Sigma $ satisfies (SM$3$%
).

Now suppose that $\Sigma _{I}\neq \emptyset $ for all nonzero ideals $I$ of
a hemiring $R$, and $B:=ker(\Sigma _{R})\neq 0$. Then, $BM=0$ for all
irreducible left $R$-semimodule $M$. Since $\Sigma _{B}\neq \emptyset $,
there exists an irreducible left $B$-semimodule $N$, and hence, $BN\neq 0$.
Whence by \cite[Lemma 4(2)]{i:otjroas} $BN$ is an irreducible left $R$%
-semimodule, and therefore, $B(BN)=0$. However, by \cite[Lemma 3]{i:otjroas} 
$BN$ is also an irreducible left $B$-semimodule, and hence, $B(BN)\neq 0$.
Thus, $ker(\Sigma _{R})=0$, and $\Sigma $ satisfies (SM$4$), too.

(ii). Using \cite[Lemma 4(1)]{i:otjroas}, the conditions (SM$3$) for $\Sigma 
$ can be established in the same fashion as in (i).

To verify (SM$4$) for the class $\Sigma $, assume that $\Sigma _{I}\neq
\emptyset $ for all nonzero ideals $I$ of a hemiring $R$ and $B:=ker(\Sigma
_{R})\neq 0$. So, $BN\neq 0$ for some a semi-irreducible left $B$-semimodule 
$N\in \Sigma _{B}$. As was shown in \cite[p. 416]{i:otjroas}, on the
semi-irreducible semimodule $N$ there exists a maximal semimodule congruence 
$\rho $ such that the $B$-semimodule $N^{^{\prime }}:=N/\rho \ $is an
reducible semimodule. Whence, $BN^{^{\prime }}\neq 0$ and, applying %
\cite[Lemma 4(2) and Corollary on p. 413]{i:otjroas}, $BN^{^{\prime }}$ is a
semi-irreducible left $R$-semimodule; and therefore, $B(BN^{^{\prime }})=0$.
However, by \cite[Lemma 4(1)]{i:otjroas} $BN^{^{\prime }}$ is also an
irreducible left $B$-semimodule, and hence, $B(BN^{^{\prime }})\neq 0$.
Thus, $ker(\Sigma _{R})=0$, and $\Sigma $ satisfies (SM$4$), too.

(iii). Suppose $ker(\Sigma _{R})=0$ for a hemiring $R$, and $I$ is a nonzero
ideal of $R$. Then, $IM\neq 0$ for some simple left $R$-semimodule $M$; and
for a nonzero element $m\in M$, the semimodule $_{R}M$ contains two
subsemimodules---$Im$ and $K:=\{x\in M$ $|$ $Ix=0\}$; and since $_{R}M$ is
simple, $K=0$ or $K=M$. However, the latter is not a case as $IM\neq 0$, and
hence, $K=0$ and $Im=M$, and therefore, the left semimodule $_{I}M$ has only
trivial $I$-subsemimodules.

Next, let $\rho $ be a congruence on $_{I}M$; and consider the congruence $%
\theta $ on $_{R}M$, where $m\theta m^{\prime }$ iff $(am,am^{\prime })\in
\rho $ for all $m,m^{\prime }\in M$ and $a\in I$. It is clear that $\theta $
is indeed a congruence on $_{R}M$ and $\rho \subseteq \theta $; and since $%
_{R}M$ is simple, $\theta =\Delta _{_{R}M}$ or $\theta =M^{2}$. If $\theta
=\Delta _{_{R}M}$, then, obviously, $\rho =\Delta _{_{I}M}$, too. If $\theta
=M^{2}$, then $(m,0)\in \theta $ and $(am,0)\in \rho $ for all $a\in I$ and
for every nonzero $m\in M$; however, as was shown above, $Im=M$, and hence, $%
(x,0)\in \rho $ for all $x\in M$. Therefore, $_{I}M$ is congruence -simple, 
\textit{i.e.}, there are only the trivial congruences on $_{I}M$, and $%
\Sigma _{I}\neq \emptyset $, and the class $\Sigma $ satisfies (SM$3$).

To verify (SM$4$) for the class $\Sigma $, consider a hemiring $R$ having $%
\Sigma _{I}\neq \emptyset $ for all nonzero ideals $I$ of $R$ and with $%
B:=ker(\Sigma _{R})\neq 0$. So, $BM=0$ for all simple left $R$-semimodules,
and there exist a simple left $B$-semimodule $X$ and an element $x\in X$
such that $Bx=X\neq 0$, and hence, $X=BX$. We may extend $_{B}X$ to a left $%
R $-semimodule $_{R}X$ defining 
\begin{equation*}
r\Sigma b_{i}x_{i}=\Sigma (rb_{i})x_{i}
\end{equation*}%
for all $b_{i}\in B,x_{i}\in X$ and $r\in R$. Indeed, if $\Sigma
b_{i}x_{i}=\Sigma b_{j}^{\prime }x_{j}^{\prime }$ ($b_{i},b_{j}^{\prime }\in
B$, $x_{i},x_{j}^{\prime }\in X$), then, since $X$ is a left $B$-semimodule, 
\begin{equation*}
b(r\Sigma b_{i}x_{i})=b\Sigma (rb_{i})x_{i}=\Sigma (brb_{i})x_{i}
\end{equation*}%
for every $b\in B$; moreover, 
\begin{equation*}
\Sigma (brb_{i})x_{i}=(br)\Sigma b_{i}x_{i}=br\Sigma b_{j}^{\prime
}x_{j}^{\prime }=\Sigma brb_{j}^{\prime }x_{j}^{\prime }=b\Sigma
rb_{j}^{\prime }x_{j}^{\prime }=b(r\Sigma b_{j}^{\prime }x_{j}^{\prime })%
\text{,}
\end{equation*}%
and so, $b(r\Sigma b_{i}x_{i})=b(r\Sigma b_{j}^{\prime }x_{j}^{\prime })$.
Consider the congruence $\alpha $ on $_{B}X$, where, for any $u,v\in X$,$\ $%
the ordered pair $(u,v)\in \alpha $ iff $bu=bv$ for all $b\in B$. Since $%
_{B}X$ is congruence-simple, $\alpha =\Delta _{_{B}X}$ or $\alpha =X^{2}$.
The latter is not a case since otherwise $(u,0)\in \alpha $ for every $u\in
X $, \textit{i.e.}, $bu=0$ for all $b\in B$, and hence, $Bx=X=0$. Hence, $%
\alpha =\Delta _{_{B}X}$. Whence, $r\Sigma b_{i}x_{i}=r\Sigma b_{j}^{\prime
}x_{j}^{\prime }$. Thus, the `scalar multiplication' by elements of $R$ is
well-defined, and it is a routine to verify that $X$ is a left $R$%
-semimodule.

Now, since $RX\supseteq BX\neq 0$, we have $RX\neq 0$, and for $X$ is a
simple left $B$-semimodule, $X$ is a simple left $R$-semimodule, too; and
hence, one gets a contradiction with $BX=0$. Therefore, there must be $%
ker(\Sigma _{R})=0$, and the class $\Sigma $ satisfies (SM$4$), too. \textit{%
\ \ \ \ \ \ }$_{\square }\medskip $

Now consider an hemiring analog of the Jacobson radical for rings introduced
by B. Bourne in \cite{b:tjroas}. First, recall \cite[Definition 3]{b:tjroas}
that a right ideal $I$ of a hemiring $R$ is \emph{right semiregular} if, for
every pair of elements $i_{1},i_{2}$ $\in I$, there exist elements $j_{1}$
and $j_{2}$ in $I$ such that $%
i_{1}+j_{1}+i_{1}j_{1}+i_{2}j_{2}=i_{2}+j_{2}+i_{1}j_{2}+i_{2}j_{1}$. By %
\cite[Theorems 3 and 4]{b:tjroas}, the sum $J(R)$ of all the right
semiregular ideals of a hemiring $R$ forms a right semiregular two-sided
ideal, which is called the \emph{Jacobson radical} of the hemiring $R$.
Also, by \cite[Theorems 5 and 6]{b:tjroas}, the mapping $\varrho :\mathbb{H}%
\longrightarrow \mathbb{H}$ given by $R\longmapsto J(R)$ is, in fact, a
radical operator in $\mathbb{H}$. Therefore, from Theorem 2.4 and %
\cite[Theorem 2.6]{hebwei:otibrtfsar}, it follows that the class 
\begin{equation*}
\mathcal{J}:=\{R\in \mathbb{H}\,|\,J(R)=R\}
\end{equation*}%
is a radical class, \emph{the Jacobson-Bourne radical class}, of $\mathbb{H}$
. Using these observations, \cite[Definitions 6 and 6', and Theorem 8]%
{i:otjroas}, Theorems 3.2 and 3.3, and Proposition 3.5, we obtain the
following hemiring analog of \cite[Example 3.14.12]{gw:rtor}:\smallskip

\noindent \textbf{Example 3.6.} \textit{Let} $\Sigma $ \textit{be one of the
classes---} $\cup _{R\in \mathbb{H}}\{_{R}M\in |_{R}\mathcal{M}|\,|\,M$ 
\textit{is irreducible}$\}$\textit{\ or }$\cup _{R\in \mathbb{H}}\{_{R}M\in
|_{R}\mathcal{M}|\,|\,M$ \textit{is semi-irreducible}$\}$--- \textit{of
semimodules over hemirings. Then,} $\Re (\Sigma )=\mathcal{J}$\textit{, the
Jacobson-Bourne radical class.\medskip }

There is another very natural hemiring analog of the Jacobson radical for
rings, namely: Applying Theorems 3.2 and 3.3 and Proposition 3.9, we obtain
that for the class of semimodules $\Sigma =\cup _{R\in \mathbb{H}}\{_{R}M\in
|_{R}\mathcal{M}|\,|\,M$ is simple$\}$, the class $\Re (\Sigma )$ is a
radical class and $\varrho _{\Re (\Sigma )}(R)=ker(\Sigma _{R})\overset{def}{%
=}$ $J_{s}(R)$ is the radical for every hemiring $R$. Obviously, on the
subclass of all rings of the class $\mathbb{H}$ both radicals, $J(R)$ and $%
J_{s}(R)$, coincide. However, in general, as the following example
demonstrates, they are different. \medskip

\noindent \textbf{Example 3.7.} Consider the Boolean semiring $\mathbf{B}%
=\{0,1\}$. It is clear that $\mathbf{B}$ is a simple $\mathbf{B}$-semimodule
and $(0:\mathbf{B})_{\mathbf{B}}=0$, and therefore, $J_{s}(\mathbf{B})=0$.
However, there are no irreducible $\mathbf{B}$-semimodules: Indeed, if $_{%
\mathbf{B}}M$ is a irreducible semimodule, then $\mathbf{B}M\neq 0$ and $%
(M,+)$ is a cancellative monoid; therefore, $1.m=(1+1)m=1.m+1.m$, and hence, 
$1.m=0$ for each $m\in M$, \textit{i.e.}, $\mathbf{B}M=0$. Thus, $J(\mathbf{B%
})=\mathbf{B}$.\medskip

As it will be shown later in Proposition 4.8, this example is a particular
case of the general observation: for commutative\ hemirings $R$, there
always takes place the inclusion $J_{s}(R)\subseteq J(R)$. In light of these
observations, it is natural to state the following, in our view interesting,
problem.\medskip

\noindent \textbf{Problem 1.} Describe the subclass of all hemirings $R$ of
the class $\mathbb{H}$ with $J_{s}(R)\subseteq J(R)$, particularly, with $%
J(R)=J_{s}(R)$.\medskip

A hemiring $R$ is called \emph{primitive }\cite[Definition 7]{i:otjroas} iff
it has a faithful irreducible left semimodule. Using Theorem 3.3 and Example
3.6, one gets the following result, removing the assumption of the additive
cancellation on primitive hemirings in \cite[Theorem 3.3, p. 12]%
{la:anotjroah} and, hence, resolving the question left open there:\medskip

\noindent \textbf{Corollary 3.8.} \textit{A hemiring} $R$ \textit{is }$J$%
\textit{-semisimple, i.e., }$J(R)=0$, \textit{iff} $R$ \textit{is
semiisomorphic to a subdirect product of primitive hemirings.}\medskip

Next, as an easy corollary of this result and our ``external'' approach to
radicals of hemirings, we immediately deduce Theorem 3.4 and Corollary 3.5
of \cite{la:anotjroah}: \medskip

\noindent \textbf{Corollary 3.9. }(\textit{cf.} \cite[Theorem 3.4, and
Corollary 3.5]{la:anotjroah}) \textit{A nonzero additively regular hemiring }%
$R$\textit{\ is }$J$\textit{-semisimple iff }$R$\textit{\ is a ring
isomorphic to a subdirect product of primitive rings.\medskip }

\noindent \textbf{Proof.} First, it is clear that any additive idempotent of
a hemiring $S$ acts on any irreducible $S$-semimodule like $0\in S$. Whence,
the zero is the only additive idempotent of a primitive additively regular
hemiring; and, therefore, any primitive additively regular hemiring is, in
fact, a ring. From this observation, Corollary 3.8, and the obvious
observation that a subdirect product of hemirings which are rings is a ring,
we end the proof. \textit{\ \ \ \ \ \ }$_{\square }\medskip $

A hemiring $R$ is \textit{congruence-simple }\cite{bshhurtjankepka:scs}
(also, \cite{bashkepka:css}, \cite{monico:ofcss}) if the diagonal, $%
\vartriangle _{R}$, and universal, $R^{2}$, congruences are the only
congruences on $R$. In the next theorem, we, in particular, give a complete
description of $J$-semisimple, congruence-simple hemirings.\medskip

\noindent \textbf{Theorem 3.10.} \textit{(i) Let }$R$\textit{\ be a
congruence-simple hemiring. Then }$R$\textit{\ is a }$J$\textit{-radical
(i.e., }$J(R)=R$) \textit{or} \textit{primitive hemiring.}

\textit{(ii) A }$J$\textit{-semisimple (i.e., }$J(R)=0$)\textit{\ hemiring }$%
R$\textit{\ is congruence-simple iff }$R$\textit{\ is a simple ring.\medskip 
}

\noindent \textbf{Proof.} (i). As is well-known (see, for example, Theorem
2.4), $J(R)$ is a subtractive ideal of $R$. But, by \cite[Proposition 4.4]%
{knt:mosssparp}, $R$ has only two trivial subtractive ideals, and therefore, 
$J(R)=R$ or $J(R)=0$. If $J(R)=0$, then by Corollary 3.8, there exists a
hemiring semiisomorphism $f$ $:R$ $\longrightarrow \prod_{i\in I}^{sub}R_{i}$%
, where $R_{i}$, $i\in I$, are primitive hemirings, and hence, $R_{i}\neq 0$%
, $i\in I$. Then, there is a nonzero surjective homomorphism $g:=\pi
_{i}|_{T}f:R\longrightarrow \prod_{i\in I}^{sub}R_{i}\longrightarrow R_{i}$
which induces the natural congruence `$\equiv $' on $R$ defined for all $%
r,s\in R$ by: $r\equiv s$ $\Longleftrightarrow $ $g(r)=g(s)$. Because $R$ is
congruence-simple, we immediately get that $g$ is injective as well, and
therefore, $R$ is isomorphic to a primitive hemiring $R_{i}$.

(ii). $\Longrightarrow $. By (i), $R$ is a primitive hemiring, and hence,
there exists a faithful irreducible left $R$-semimodule $M$. By %
\cite[Proposition 3.1]{zumbr:cofcsswz}, congruence-simpleness of a hemiring $%
R$\ implies that either $R$ is a simple ring or the reduct $(R,+,0)$ is an
idempotent monoid. In the latter case, because $(M,+)$ is a cancellative
monoid, $RM=0$, and hence, the semimodule $M$ is not irreducible. Thus, $R$
is a simple ring.

$\Longleftarrow $. It is obvious. \textit{\ \ \ \ \ \ }$_{\square }\medskip $

Let $M$ be a join semilattice with zero; and following \cite{zumbr:cofcsswz}%
, for any elements $a,b\in M$, consider the following mappings on $M$: 
\begin{equation*}
e_{a,b}(x):=\{_{b\text{ otherwise \ }}^{0\text{ if }x\leq a,}(x\in M)\text{.}
\end{equation*}%
As was shown in \cite[Lemma 2.2]{zumbr:cofcsswz}, $e_{a,b}\in End(M)$ for
all $a,b\in M$. Then, let $\mathbf{F}_{M}$ be the submonoid of $(End(M),+,0)$
generated by the endomorphisms $e_{a,b}$, $a,b\in M$. It is easy to see
(see, also, \cite[Lemma 2.2]{zumbr:cofcsswz} and \cite[Lemma 3.1]%
{knz:ososacs}) that $\mathbf{F}_{M}$ is a left ideal of $End(M)$, which is
an ideal when $M$ is finite. A subhemiring $R$ of $End(M)$ is said to be 
\emph{dense} \cite{zumbr:cofcsswz} (also, \cite[p. 154]{irs:rtofsos}) if it
contains $\mathbf{F}_{M}$. We conclude this subsection with a
characterization of finite additively-idempotent $J_{s}$-semisimple
hemirings. \medskip

\noindent \textbf{Theorem 3.11.} \textit{A finite additively idempotent
hemiring} $R$ \textit{is} $J_{s}$\textit{-semisimple (i.e., } $J_{s}(R)=0$%
\textit{)} \textit{iff} \textit{it} \textit{is semiisomorphic to a subdirect
product of some hemirings} $S_{i}(i\in I)$ \textit{such that each} \textit{%
of the hemirings} $S_{i}(i\in I)$\textit{, in turn, is isomorphic to a dense
subhemiring of} \textit{the endomorphism hemiring} $End(M_{i})$ $(i\in I)$ 
\textit{of a finite semilattice} $M_{i}$ $(i\in I)$ \textit{with zero.}%
\medskip

\noindent \textbf{Proof.} $\Longleftarrow $. By \cite[Corollary 5.4]%
{irs:rtofsos}, for each $i\in I$, there is a simple left $S_{i}$-semimodule $%
M_{i}$ such that the map $\varphi :S_{i}\longrightarrow End(M_{i},+)$,
defined for all $s\in S_{i}$ and $m\in M_{i}$ by $\varphi (s)(m)=sm$, is an
injective homomorphism of hemirings, and hence, $(0:M_{i})_{S_{i}}=0$.
Whence by Theorem 3.3, $J_{s}(R)=0$.

$\Longrightarrow $. By Theorem 3.3, the hemiring $R$ is semiisomorphic to a
subdirect product of finite additively-idempotent hemirings $R_{i}(i\in I)$
such that, for each $i\in I$, there exists a simple faithful left $R_{i}$%
-semimodule $M_{i}$ for which, of course, $(M_{i},+)$ is a finite idempotent
monoid and $M_{i}=R_{i}m$ for any nonzero $m\in M_{i}$. Thus, there exist
the hemiring homomorphism $\varphi _{i}:R_{i}\longrightarrow End(M_{i},+)$
defined for all $r\in R_{i}$ and $m\in M_{i}$ by $\varphi (r)(m)=rm$ and
having $ker(\varphi _{i}):=\varphi _{i}^{-1}(0)=(0:M_{i})_{R_{i}}=0$, as
well as the corresponding congruence ``$\equiv _{\varphi _{i}}$'' on $R_{i}$
such that $r\equiv _{\varphi _{i}}r^{\prime }$ $\Longleftrightarrow $ $%
\varphi _{i}(r)=\varphi _{i}(r^{\prime })$ for all $r,r^{\prime }\in R_{i}$
and, hence, the injective homomorphism $\psi :R_{i}/\equiv _{\varphi
_{i}}\longrightarrow End(M_{i},+)$ such that $\psi (\overline{r})=\varphi
_{i}(r)$ for all $r\in R_{i}$. For the natural projection $\pi
:R_{i}\longrightarrow R_{i}/\equiv _{\varphi _{i}}$, obviously, $ker(\pi
):=\pi ^{-1}(\overline{0})=(0:M_{i})_{R_{i}}=0$, and therefore, $\pi $ is a
semiisomorphism and $R$ is semiisomorphic to a subdirect product of the
hemirings $S_{i}:=R_{i}/\equiv _{\varphi _{i}}$. Also, it is clear that $%
M_{i}$ with the `scalar multiplication' given by $\overline{r}m=rm$ for all $%
m\in M$ and $r\in R_{i}$ is a finite simple left $S_{i}$-semimodule. Now,
for the homomorphism $\psi :R_{i}/\equiv _{\varphi _{i}}\longrightarrow
End(M_{i},+)$ is injective and \cite[Theorem 5.3]{irs:rtofsos}, we conclude
that $S_{i}$ is isomorphic to a dense subhemiring of the endomorphism
hemiring $End(M_{i})$ of a finite semilattice $M_{i}$ with zero. \textit{\ \
\ \ \ \ }$_{\square }\medskip $

\subsection{The Brown-McCoy Radical of Semirings}

The Brown-McCoy radical for hemirings constitutes another important

\noindent hemiring analog of a classical radical for rings. Consider the
subclass of ideal-simple semirings $\mathbb{S}:=\{R\in \mathbb{H}\,|\,R$ is
an ideal-simple semiring$\}$ of the class $\mathbb{H}$. It is obvious that
subclass $\mathbb{S}$ is a regular class of $\mathbb{H}$, and therefore, by
Theorem 2.3, the class $\mathcal{U}\mathbb{S}=\{R\in \mathbb{H}\,|\,R$ has
no nonzero homomorphic images in $\mathbb{S}\}$ forms a radical class, the 
\emph{Brown-McCoy radical class}, of $\mathbb{H}$; and\ by $\mathcal{R}%
_{BM}(R)$ we will denote the corresponding Brown-McCoy radical of a hemiring 
$R$. Because for any semiring $R$, obviously, there exists a maximal
congruence $\rho $ on $R$ such that $R/\rho $ is a congruence-simple
semiring, it is clear that a hemiring $R$ is a Brown-McCoy radical hemiring, 
\textit{i.e.} $\mathcal{R}_{BM}(R)=R$, iff there is no nonzero homomorphic
images of $R$ in the class of all simple semirings.

The Brown-McCoy radicals for rings and some other algebraic systems can be
characterized in terms of some special congruences (see, for instance, %
\cite[Section 4.8]{gw:rtor} and \cite[p. 430]{la:mcatbrfs} and references
there). In our next result, we present such a characterization for hemirings
belonging to a fairly large class --- including, in particular, all
commutative hemirings and rings --- of hemirings. For that we need natural
hemiring analogs of some well-known for rings notions.

Given a hemiring $R$ and element $r\in R$, consider two congruences, $\rho
_{G(r)}$ and $\rho _{F_{1}(r)}$, on $R$, \ generated by two subsets 
\begin{eqnarray*}
G(r) &:&=\{(rx+\Sigma _{i=1}^{n}y_{i}rz_{i},\text{ }x+\Sigma
_{i=1}^{n}y_{i}z_{i})\,|\,n\in \mathbb{N},\text{ }x,y_{i},z_{i}\in R\}\text{
and} \\
F_{1}(r) &:&=\{(rx+yr,\text{ }x+y)\,|\,x,y\in R\}
\end{eqnarray*}%
of $R^{2}$, respectively. Further, following \cite[Section 4.8]{gw:rtor}, we
say that the element $r$ is $G$\emph{-regular} ($F_{1}$\emph{-regular}) if $%
(r,0)\in \rho _{G(r)}$ ($(r,0)\in \rho _{F_{1}(r)}$); and a hemiring $R$ is
called $G$\emph{-regular} ($F_{1}$\emph{-regular}) if every element $r\in R$
is $G$\emph{-regular} ($F_{1}$\emph{-regular}). Now we are ready to extend
the characterization of the Brown-McCoy radical for rings \cite[Theorem
4.8.2]{gw:rtor} to the subclasses of subtractive, commutative, and
lattice-ordered, hemirings.\medskip

\noindent \textbf{Theorem 3.12.} \textit{Let} $R$ \textit{be a subtractive
(or commutative, or lattice-ordered) hemiring. Then, the following
statements are equivalent:}

\textit{(i)} $R$ \textit{is a Brown-McCoy radical hemiring, i.e.,} $R\in 
\mathcal{U}\mathbb{S}$;

\textit{(ii)} $R$ \textit{is} $G$\textit{-regular;}

\textit{(iii)} $R$ \textit{is} $F_1$\textit{-regular.\medskip}

\noindent \textbf{Proof.} (i) $\Longrightarrow $ (ii). Let $(r,0)\notin $ $%
\rho _{G(r)}$ for some element $r\in R$. By the Zorn's Lemma, there exists a
maximal congruence $\rho $ on $R$ such that $\rho _{G(r)}\subseteq \rho $
and $(r,0)\notin $ $\rho $. In fact, $\rho $ is a maximal congruence on $R$:
Indeed, since $(r,0)\in \theta $ for any properly containing $\rho $
congruence $\theta $ on $R$ and $(rx,0),(rx,x)\in \theta $ for all $x\in R$,
one has $(x,0)\in \theta $ for all $x\in R$, and hence, $R^{2}\subseteq
\theta $. So, $S=R/\rho $ is a congruence-simple hemiring with $(r,0)\notin
\rho $, \textit{i.e.}, $\overline{r}\neq \overline{0}$. For $(rx,x)\in
G(r)\subseteq \rho $, we have $\overline{r}\,\overline{x}=\overline{x}$ for
all $x\in R$, and hence, $\overline{r}$ is a left identity on $S$. We claim
that, actually, $\overline{r}$ is an identity on $S$. Indeed, for the
relation $\gamma :=\{(x,y)\in S^{2}\,|\,x\overline{r}=y\overline{r}\}$ on $S$%
, it is easy to see that $(x+a,y+a)\in \gamma $ and $(ax,ay)\in \gamma $ for
any $(x,y)\in \gamma $ and $a\in S$; moreover, since $x\overline{r}=y%
\overline{r}$ and $\overline{r}a=a$, we have $xa=x(\overline{r}a)=(x%
\overline{r})a=(y\overline{r})a=y(\overline{r}a)=ya$, and therefore, $%
(xa,ya)\in \gamma $ and, hence, $\gamma $ is a congruence on $S$. Thus, $%
\gamma =\vartriangle _{S}$ or $\gamma =S^{2}$. In the latter, however, $(%
\overline{r},\overline{0})\in \gamma $ and, hence, $\overline{r}=\overline{r}%
\,\overline{r}=\overline{0}\overline{r}=\overline{0}$ what contradicts with $%
\overline{r}\neq \overline{0}$. Whence, $\gamma =\vartriangle _{S}$. As
obviously $(y\overline{r},y)\in \gamma =\vartriangle _{S}$ for all $y\in S$,
one gets that $y\overline{r}=y$ for all $y\in S$ and $\overline{r}$ is also
a right identity on $S$. Now consider the following three cases.

If $R$ is a commutative hemiring, then $S$ is a commutative
congruence-simple semiring, and therefore, by \cite[Theorem 3.2]{mf:ccs} (or %
\cite[Theorem 10.1]{bshhurtjankepka:scs}), $S$ is either a field or
isomorphic to the Boolean semifield $\mathbf{B}:=\{0,1\}$. Anyway, $S$ is an
ideal-simple semiring.

If the hemiring $R$ is subtractive hemiring, then $S$ is also subtractive,
and hence, a subtractive congruence-simple semiring, which, by %
\cite[Proposition 4.4]{knt:mosssparp}, is ideal-simple, too.

Finally, let $R$ be a lattice-ordered hemiring. Then, $(R,+,0)$ and $(S,+,%
\overline{0})$ are additively idempotent monoids; and let $I\subseteq S$ be
a nonzero ideal of $S$. For $S$ is congruence-simple, the Bourne congruence
``$\equiv _{I}$'' on $S$ coincides with $S^{2}$. So, $\overline{r}\equiv _{I}%
\overline{0}$, \textit{i.e.}, there exist elements $a,b\in R$ such that $%
\overline{a},\overline{b}\in I$ and $\overline{b}=$ $\overline{r}+$ $%
\overline{a}$. For $\overline{a}=\overline{r}\overline{a}=\overline{ra}\leq 
\overline{r\wedge a}\leq \overline{r}$, it follows that $\overline{b}=$ $%
\overline{r}+$ $\overline{a}=\overline{r}\in I$, and hence, $I=S$; and
therefore, $S$ is an ideal-simple semiring in this case, too.

Thus in all three cases, we have shown that $R$ has a nonzero homomorphic
image $S$ which is a simple semiring, and therefore, $R\notin \mathcal{U}%
\mathbb{S}$.

(ii) $\Longrightarrow $ (iii). It is obvious since $\rho _{G(r)}$ $\subseteq 
$ $\rho _{F_{1}(r)}$ for all $r\in R$.

(iii) $\Longrightarrow $ (i). In fact, this implication is true for any
hemiring $R$, because if $R$ is a $F_{1}$-regular hemiring and $\tau $ is a
congruence on $R$ such that $\tau \neq R^{2}$, then the factor hemiring $%
S=R/\tau $ has no identity element: Indeed, if for some $e\in R$ the element 
$\overline{e}\in S$ is an identity of $S$, then, as clearly $%
F_{1}(e)\subseteq \tau $, one has that $(e,0)\in \rho _{F_{1}(e)}\subseteq
\tau $. Thus, $R$ has no nonzero homomorphic images in $\mathbb{S}$, and
therefore, $\mathcal{R}_{BM}(R)=R$. \textit{\ \ \ \ \ \ }$_{\square
}\medskip $

In contrast to rings, the implication (i) $\Longrightarrow $ (iii) in
Theorem 3.12 in general is not true for hemirings, namely: \medskip

\noindent \textbf{Example 3.13.} \textit{A semiring }$End(M)$\textit{\ of
endomorphisms of a}

\noindent \textit{non-distributive finite lattice }$M$\textit{\ is not }$%
F_{1}$\textit{-regular semiring with}

\noindent $\mathcal{R}_{BM}(End(M))=End(M)$\textit{.\medskip }

\noindent \textbf{Proof.} By \cite[Proposition 2.3]{zumbr:cofcsswz}, $End(M)$
is a congruence-simple semiring, and therefore, any nonzero surjective
homomorphism $\alpha :End(M)\longrightarrow S$ to an ideal-simple semiring $%
S $, would be an isomorphism, and hence, $End(M)$ would be a simple
semiring; the latter, in the notations introduced before Theorem 3.11, would
imply that $\mathbf{F}_{M}=End(M)$, what, by \cite[Proposition 4.9 and
Remark 4.10]{zumbr:cofcsswz}, would be a contradiction with that $M$ is a
non-distributive lattice. Thus, there is no nonzero surjective homomorphism $%
\alpha :End(M)\longrightarrow S$ to an ideal-simple semiring $S$, and
therefore, $\mathcal{R}_{BM}(End(M))=End(M)$.

However, $F_{1}(id_{M})=\{(id_{M}f+gid_{M},$ $f+g)\,|\,$\ $f,g\in End(M)\}$ $%
=\vartriangle _{End(M)}$. So, the congruence $\rho _{F_{1}(id_{M})}=$ $%
\vartriangle _{End(M)}$, and therefore, $(id_{M},0)\notin \rho
_{F_{1}(id_{M})}$ and $End(M)$ is not $F_{1}$-regular semiring.\textit{\ \ \
\ \ \ }$_{\square }\medskip $

In light of these observations, it seems to be natural to bring up the
following problem.\medskip

\noindent \textbf{Problem 2.} Describe all hemirings for which the
implication (i) $\Longrightarrow $ (ii) in Theorem 3.12 is true.\medskip

Extending the corresponding result for rings (see, \textit{i.g.}, %
\cite[Theore. 4.8.1]{gw:rtor}), in the next theorem we present an
alternative, ``working,'' description of the Brown-McCoy radical of
hemirings from a fairly large class --- including, in particular, all
commutative hemirings and rings --- of hemirings. To do that, we first note
the following, generally useful and more or less obvious, fact.\medskip

\noindent \textbf{Lemma 3.14.} \textit{If }$S$\textit{\ is an ideal of a
hemiring }$R$\textit{\ and }$\rho $\textit{\ is a congruence on }$S$\textit{%
\ such that }$S/\rho $\textit{\ is a semiring with the identity element }$%
\overline{e}$\textit{, then the relation }$\overline{\rho }$\textit{\ on }$R$%
\textit{, where }$(a,b)\in \overline{\rho }\ $iff $(eae,ebe)\in \rho $%
\textit{\ for all }$a,b\in R$,\textit{\ is a congruence on }$R$\textit{\ and 
}$\rho \subseteq \overline{\rho }$\textit{. Furthermore, the natural map }$%
\varphi :R/\overline{\rho }\longrightarrow S/\rho $\textit{\ given by }$%
\overline{r}\longmapsto \overline{ere}$\textit{\ is a hemiring isomorphism;
in particular, if }$R$\textit{\ is a semiring, then }$\varphi $\textit{\ is
a semiring isomorphism.\medskip }

\noindent \textbf{Proof.} It is clear that $\overline{\rho }$ is an
equivalence relation on $R$. Now, let $(a,b)\in \overline{\rho }$ and $c\in
R $. For $e\in S$ and $S$ is an ideal of $R$, one has $(eae,ebe)\in \rho $
and $ece\in S$, and hence, $(e(a+c)e,e(b+c)e)=(eae,ebe)+(ece,ece)\in \rho $;
and therefore, $(a+c,b+c)\in \overline{\rho }$. Also, from $ec,ae\in S$ and $%
(eae,ebe)\in \rho $, we have $(eceae,ecebe)\in \rho $. Since $\overline{e}$
is an identity of $S/\rho $, one has $(ece,ec)\in \rho $, and hence, $%
(eceae,ecae)\in \rho $. Similarly, one gets $(ecebe,ecbe)\in \rho $. Whence, 
$(ecae,ecbe)\in \rho $ and $(ca,cb)\in \overline{\rho }$. In the similar
way, one can show that $(ac,bc)\in \overline{\rho }$, too. Thus, $\overline{%
\rho }$ is a congruence on $R$. The rest is trivial.\textit{\ \ \ \ \ \ }$%
_{\square }\medskip $

For any hemiring $R$, let $\mathfrak{C}_{s}$ stay for the set of all
congruence $\rho $ on $R$ such that $R/\rho $ is a nonzero simple semiring,
and $\rho _{s}:=$ $\cap \{\rho \,|\,\rho \in \mathfrak{C}_{s}\}$. For each $%
\rho \in \mathfrak{C}_{s}$, by $[0_{\rho }]$ we denote the \textit{kernel of 
}$\rho $, \textit{i.e., }$[0_{\rho }]:=\{r\in R$ $|$ $(r,0)\in \rho
\}\subseteq R$; and let $\mathcal{R}_{\mathfrak{C}_{s}}(R):=\cap \{[0_{\rho
}]\,|\,\rho \in \mathfrak{C}_{s}\}$. It is easy to see that $\mathcal{R}_{%
\mathfrak{C}_{s}}(R)$ is a subtractive ideal of $R$, and $a\in \mathcal{R}_{%
\mathfrak{C}_{s}}(R)$ iff $(a,0)\in \rho _{s}$, \textit{i.e.}, $\mathcal{R}_{%
\mathfrak{C}_{s}}(R)$ is the kernel of $\rho _{s}$.

Also, we say that a hemiring $R$ is \emph{strongly subtractive} if every
ideal $I$ of $R$ is a subtractive hemiring. The class of strongly
subtractive hemirings, obviously, includes all rings, but not only them,
namely (\textit{cf}. \cite[Fact 2.1]{knt:mosssparp}): Let $R$ be the chain $%
0<a<b<1$ with the operation of multiplication defined as follows: $%
b^{2}=b,\,a^{2}=0,\,ba=a,\,ab=0$; then, $(R,\vee ,\cdot )$ is a semiring and 
$\{0\},\{0,a\},\{0,a,b\}\ $and $R$ are the only ideals of $R$, and they all
are clearly subtractive hemirings.\medskip

\noindent \textbf{Theorem 3.15.} \textit{Let} $R$ \textit{be a commutative
(or strongly subtractive, or lattice-ordered) hemiring. Then, }$\mathcal{R}%
_{BM}(R)=\mathcal{R}_{\mathfrak{C}_{s}}(R)$\textit{.}$\medskip $

\noindent \textbf{Proof.} For each $[0_{\rho }]$, $\rho \in \mathfrak{C}_{s}$%
, is a subtractive ideal of $R$ that is the kernel of the natural projection 
$\pi :R\longrightarrow R/\rho $ and \cite[Proposition 10.16]{golan:sata}, $%
R/[0_{\rho }]$ and $R/\rho $ are semiisomorphic. As by \cite[Theorem 2.6(c)]%
{hebwei:otibrtfsar} $\mathbb{S}=\{R\in \mathbb{H}\,|\,\mathcal{R}%
_{BM}(R)=0\} $ is a semisimple class of $\mathbb{H}$, and $\mathcal{R}%
_{BM}(R/\rho )=0$, \textit{i.e.}, $R/\rho \in \mathbb{S}$ since $R/\rho $ is
a simple semiring, from \cite[Theorem 3.7(b)]{hebwei:scos} it follows that $%
R/[0_{\rho }]\in \mathbb{S}$, \textit{i.e.}, $\mathcal{R}_{BM}(R/0_{\rho
})=0 $, too. Whence, taking into consideration that by \cite[Theorem 4.9]%
{m:otrtfs} $\mathcal{R}_{BM}(R)=\cap \{K\in \mathcal{SI}(R)\,|\,\mathcal{R}%
_{BM}(R/K)=0\}$, we have $\mathcal{R}_{BM}(R)\subseteq \cap \{0_{\rho
}\,|\,\rho \in \mathfrak{C}_{s}\}=\mathcal{R}_{\mathfrak{C}_{s}}(R)$.

Now let $S$ stay for $\mathcal{R}_{\mathfrak{C}_{s}}(R)$. Using the
A-D-S-property of a radical class (see, \cite[Theorem 6.2]{m:otrtfs}), we
have $\mathcal{R}_{BM}(S)\subseteq \mathcal{R}_{BM}(R)$ for the ideal $%
\mathcal{R}_{BM}(S)$ of $R$. Suppose that $\mathcal{R}_{BM}(S)\neq S$. It is
clear that $S$ is a commutative, or subtractive, hemiring as soon as $R$ is
a commutative, or strongly subtractive, hemiring itself. Also, if $R$ is a
lattice-ordered hemiring, $S$ is a lattice-ordered hemiring, too: Indeed,
for any $a,b\in S$, we have $a+(a\wedge b)=a\vee (a\wedge b)=a$, and as $S$
is a subtractive ideal of $R$ and $a\in S$, we conclude that $a\wedge b\in S$%
. From these observations and Theorem 3.12, it follows that there exists an
element $e\in S$ such that $e$ $\notin $ $\rho _{F_{1}(e)}$. For the Zorn's
Lemma, let $\rho $ be a maximal congruence on $S$ such that $%
F_{1}(e)\subseteq \rho $ and $(e,0)\notin \rho $. As in Theorem 3.12, it can
be shown that $S/\rho $ is a simple semiring with the identity $\overline{e}$
$\in $ $S/\rho $ for $e\in S$. Applying Lemma 3.14, we have the congruence $%
\overline{\rho }$ on $R$ such that for all $a,b\in R$, 
\begin{equation*}
(a,b)\in \overline{\rho }\Longleftrightarrow (eae,ebe)\in \rho \text{,}
\end{equation*}%
and $\rho \subseteq \overline{\rho }$. It is clear that $(e,0)\notin 
\overline{\rho }$ and $(ex,x),\,(xe,x)\in \overline{\rho }$ for all $x\in R$%
; and hence, $F_{1}(e):=\{(ex+ye,x+y)\,|\,x,y\in R\}\subseteq \overline{\rho 
}$. Again, for the Zorn's Lemma, let $\delta $ be a maximal congruence on $R$
such that $F_{1}(e)\subseteq \delta $ and $(e,0)\notin \delta $; and as in
Theorem 3.12, we get that $R/\delta $ is a simple semiring with $(e,0)\notin
\delta $ that is a contradiction with $e\in \mathcal{R}_{\mathfrak{C}%
_{s}}(R)=S$. Therefore, $\mathcal{R}_{BM}(S)=S$ and $\mathcal{R}_{\mathfrak{C%
}_{s}}(R)\subseteq \mathcal{R}_{BM}(R)$.\textit{\ \ \ \ \ \ }$_{\square
}\medskip $

Applying Theorem 3.15, we obtain the following complete descriptions of $%
\mathcal{R}_{BM}$-semisimple commutative, or lattice-ordered,
hemirings.\medskip

\noindent \textbf{Corollary 3.16.} \textit{For a commutative hemiring} $R$%
\textit{,} \textit{the following conditions are equivalent:}

\textit{(i )} $\mathcal{R}_{BM}(R)=0$\textit{;}

\textit{(ii)} $\mathcal{R}_{\mathfrak{C}_{s}}(R)=0$\textit{;}

\textit{(iii)} $R$ \textit{is semiisomorphic to a subdirect product of
simple commutative semirings;}

\textit{(iv)} $R$ \textit{is semiisomorphic to a subdirect product of
semirings that are either fields or the semifield} $\mathbf{B}$\textit{%
.\medskip }

\noindent \textbf{Proof.} (i) $\Longleftrightarrow $ (ii) by Theorem 3.15.

The equivalence (iii) $\Longleftrightarrow $ (iv) follows from the fact (%
\cite[Theorem 3.2]{mf:ccs}, or \cite[Thorem 10.1]{bshhurtjankepka:scs}) that
a simple commutative semiring is either a field or the semifield $\mathbf{B}$%
.

(ii) $\Longrightarrow $ (iii). As $\mathcal{R}_{\mathfrak{C}_{s}}(R)=0$, one
has that $R$ is semiisomorphic to a subdirect product $\prod_{\rho \in 
\mathfrak{C}_{s}}^{sub}R/\rho $ of simple commutative semirings\ $R/\rho $, $%
\rho \in \mathfrak{C}_{s}$.

(iii) $\Longrightarrow $ (i). As simple semirings are obviously $\mathcal{R}%
_{BM}$-semisimple, the implication immediately follows from \cite[Theorem
4.3 and Theorem 3.7 (b)]{hebwei:scos}.\textit{\ \ \ \ \ \ }$_{\square
}\medskip $

\noindent \textbf{Corollary 3.17.} \textit{For a lattice-ordered hemiring} $%
R $\textit{,} \textit{the following conditions are equivalent:}

\textit{(i )} $\mathcal{R}_{BM}(R)=0$\textit{;}

\textit{(ii)} $\mathcal{R}_{\mathfrak{C}_{s}}(R)=0$\textit{;}

\textit{(iii)} $R$ \textit{is semiisomorphic to a subdirect product of
simple lattice-ordered semirings;}

\textit{(iv)} $R$ \textit{is semiisomorphic to a subdirect product of the
Boolean semifields} $\mathbf{B}. \medskip$

\noindent \textbf{Proof.} Using the fact (see, \cite[Theorem 6.7]%
{knz:ososacs}) that a lattice-ordered semiring $S$ is simple iff $S\cong 
\mathbf{B}$, this result can be proved in the similar fashion as in
Corollary 3.16. \textit{\ \ \ \ \ \ }$_{\square }\medskip $

Let $End(M)$ be an endomorphism semiring of a non-distributive finite
lattice $M$. By \cite[Proposition 2.3]{zumbr:cofcsswz}, this
congruence-simple semiring is not commutative, not subtractive, and not
lattice-ordered. Also, as was shown in Example 3.13, it has no nonzero
surjective homomorphisms onto simple semirings, and therefore, $\mathcal{R}_{%
\mathfrak{C}_{s}}(End(M))=End(M)=$ $\mathcal{R}_{BM}(End(M))$. In light of
this observation and Theorem 3.15, we conclude this section with the
following, as we think interesting, question: \medskip

\noindent \textbf{Problem 3.} Is it true that $\mathcal{R}_{BM}(R)=\mathcal{R%
}_{\mathfrak{C}_{s}}(R)$ for all hemirings $R\in \mathbb{H}$? If the answer
is ``NO,'' then describe all hemirings $R\in \mathbb{H}$ for which $\mathcal{%
R}_{BM}(R)=\mathcal{R}_{\mathfrak{C}_{s}}(R)$.

\section{The Nakayama's Lemma and Jacobson - \ \ \ \ \ \ \ \ \ \ \ \ \ \ \ \
\ \ \ \ \ \ \ \ \ \ \ \ \ \ \ \ \ \ \ \ \ \ \ \ \ \ \ Chevalley Density
Theorem for Semirings}

In this section, we establish semiring analogs of the well-known classical
ring results---Nakayama's and Hopkins Lemmas and Jacobson-Chevalley Density
Theorem---by reducing our semiring settings to corresponding\ original ring
ones. But to do that and for the reader's convenience, we first recall some
notions we need in a sequel.

On any left semimodule $M$ over a hemiring $R$, there exists the congruence,
``$\equiv $'', defined for all $m,m^{\prime }\in M$ as follows: $m\equiv
m^{\prime }$ iff $m+x=m^{\prime }+x$ for some $x\in M$. Let $M^{\ast }$ and $%
m^{\ast }\in $ $M^{\ast }$ stay for the factor semimodule $M/\equiv $ and
equivalence class of an element $m\in M$, respectively. In particular,
considering a hemiring $R$ as the regular semimodule $_{R}R$, it is easy to
see that $R^{\ast }$ becomes a hemiring with the multiplication $r^{\ast
}s^{\ast }=(rs)^{\ast }$ for all $r^{\ast },s^{\ast }\in R^{\ast }$. It is
obvious that $_{R}M\cong $ $_{R}M^{\ast }$ for any cancellative semimodule $%
M $, and $M^{\ast }$ is also a cancellative left $R^{\ast }$-semimodule with
the scalar multiplication defined for all $m^{\ast }\in M^{\ast }$ and $%
r^{\ast }\in R^{\ast }$ by $r^{\ast }m^{\ast }=(rm)^{\ast }$. Furthermore,
it is easy to see that $M^{\ast }$ is an irreducible left $R^{\ast }$%
-semimodule for any irreducible left $R$-semimodule $M$, and if $M$ is an
irreducible left $R^{\ast }$-semimodule, then $M$ is also an irreducible
left $R$-semimodule with the scalar multiplication given by $rm=r^{\ast }m$
for all $r\in R$ and $m\in M$ . From these remarks and \cite[Thorem 8]%
{i:otjroas}, it follows that $J(R)=\{r\in R\,|\,r^{\ast }\in J(R^{\ast })\}$
(see, also \cite[p. 420]{i:otjroas}).

For any left $R$-semimodule $M$, there exists the left $R$-\emph{module of
differences} $D(M)$ of $M$ \cite[Chapter 16]{golan:sata} defined as the
factor semimodule of the left $R$-semimodule $M\times M$ with respect to the
subsemimodule $W=\{(m,m)$ $|$ $m\in M$ $\}\subseteq $ $M\times M$, \textit{%
i.e.}, $D(M):=(M\times M)/W$. In fact, the semimodule $D(M)$ is a left $R$%
-module since for any $(m,m^{\prime })\in M\times M$ in $D(M)$ one has $%
\overline{(m,m^{\prime })}+\overline{(m^{\prime },m)}=\overline{(0,0)}$.
Also, there exists the canonical $R$-homomorphism $\xi _{M}:$ $M$ $%
\longrightarrow D(M)$ given by $m\longmapsto \overline{(m,0)}$. In the case
when $M$ is a cancellative semimodule, $\xi _{M}$ is injective, and
therefore, we can consider the elements $\overline{(m,0)}$ and $m$ to be the
same and any element $\overline{(m,m^{\prime })}\in $ $D(M)$ to be the
``difference'' of the elements $\overline{(m,0)}$ and $\overline{%
(0,m^{\prime })}$, \textit{i.e.}, $D(M)=\{m-m^{\prime }\,|\,m,m^{\prime }\in
M\}$. In particular, the left $R$-module of differences $D(R)$ of the
regular semimodule $_{R}R$ can be considered as a ring --- the \emph{ring of
differences} of $R$ \cite[Chapter 8, p. 101]{golan:sata} --- with the
operation of multiplication defined for all $a,b,c,d\in R$ by $\overline{%
(a,b)}\overline{(c,d)}=\overline{(ac+bd,ad+cb)}$; and if $R$ is a semiring,
then the ring of differences $D(R)$ is also a semiring with the identity $%
\overline{(1,0)}$. Moreover, it is easy to see that $D(M)$ becomes a left $%
D(R)$-module with $\overline{(a,b)}$ $\overline{(m_{1},m_{2})}=\overline{%
(am_{1}+bm_{2},am_{2}+bm_{1})}$ for all $a,b\in R$ and $m_{1},m_{2}\in M$,
and $_{R}D(M)\cong $ $_{R}D(M^{\ast })$, as well as $D(R)\cong D(R^{\ast })$
as rings. Then, it is easy to see (see also \cite[p. 419, Section 4 c)]%
{i:otjroas}) that a cancellative left $R$-semimodule $M$ is irreducible iff $%
D(M)$ is an irreducible left $D(R)$-module, as well as $J(R^{\ast
})=J(D(R^{\ast }))\cap R^{\ast }=J(D(R))\cap R^{\ast }$ for any hemiring $R$ %
\cite[p. 420, Section 4 e)]{i:otjroas}.

Now we introduce the concept of a weakly-finitely (in short, $w$-finitely)
generated semimodule, which is a very natural\ generalization of the notion
of an irreducible semimodule and will prove to be useful in the ``reduction
procedure'' from our semiring settings to corresponding ring ones in a
sequel.\medskip

\noindent \textbf{Definition 4.1.} A left $R$-semimodule $M$ over a hemiring 
$R$ is called $w$\emph{-finitely generated} if there exist a natural number $%
n$ and pairs $(m_{1},m_{1}^{\prime }),$ $(m_{2},m_{2}^{\prime
}),...,(m_{n},m_{n}^{\prime })\in $ $M^{2}$ such that, for any element $m\in
M$, there exist pairs $(r_{1},s_{1}),(r_{2},s_{2}),...,(r_{n},s_{n})\in
R^{2} $ and 
\begin{equation*}
m+(r_{1}m_{1}^{\prime }+s_{1}m_{1})+...+(r_{n}m_{n}^{\prime
}+s_{n}m_{n})=(r_{1}m_{1}+s_{1}m_{1}^{\prime
})+...+(r_{n}m_{n}+s_{n}m_{n}^{\prime })\text{.}
\end{equation*}

\noindent \textbf{Proposition 4.2.} \textit{(i) Any irreducible left }$R$%
\textit{-semimodule is }$w$\textit{-finitely generated as well.}

\textit{(ii) Any finitely generated left }$R$\textit{-semimodule is }$w$%
\textit{-finitely generated as well.}

\textit{(iii) For any }$w$\textit{-finitely\ generated left }$R$\textit{%
-semimodule\ }$M$\textit{, the semimodule }$M^{\ast }$\textit{\ is an }$w$%
\textit{-finitely generated left }$R^{\ast }$\textit{-semimodule.}

\textit{(iv) For any }$w$\textit{-finitely\ generated left }$R$\textit{%
-semimodule\ }$M$\textit{,\ the left }$D(R)$\textit{-module }$D(M)$\textit{\
is finitely generated.}

\textit{(v) A cancellative left }$R$\textit{-semimodule }$M$\textit{\ is }$w$%
\textit{-finitely generated iff the left }$D(R)$\textit{-module }$D(M)$ 
\textit{is finitely generated.\medskip }

\noindent \textbf{Proof.} (i), (ii) and (iii) are clear.

(iv). Let $M$ be an $w$-finitely generated left $R$-semimodule, and

\noindent $(m_{1},m_{1}^{\prime }),(m_{2},m_{2}^{\prime
}),...(m_{n},m_{n}^{\prime })$ $\in M^{2}$ be the pairs as in Definition
4.1. Then one may easily see that $D(M)=D(R)\overline{(m_{1},m_{1}^{\prime })%
}+...+D(R)\overline{(m_{n},m_{n}^{\prime })}$, and therefore, $D(M)$ is a
finitely generated left $D(R)$-module.

(v). This follows from (iv) and the fact that $D(M)=\{m-m^{\prime
}\,|\,m,m^{\prime }\in M\}$ when $M$ is cancellative. \textit{\ \ \ \ \ \ }$%
_{\square }\medskip $

As usually (see, \textit{e.g.}, \cite[pages 50, 155]{golan:sata}), $Z(R)=\{$ 
$r\in R\,|\,r+x=x$ for some $x\in R$ $\}$ and $Z(M)=\{$ $m\in M\,|\,m+x=x$
for some $x\in M$ $\}$ denote the \emph{zeroids} of a hemiring $R$ and an $R$%
-semimodule $M$, respectively. We are now ready to present a semiring
version of the fundamental in the theory of rings and modules, famous
Nakayama's Lemma (see, for instance, \cite[Nakayama's Lemma 4.22]{lam:afcinr}%
).\medskip

\noindent \textbf{Theorem 4.3.} \textit{Let }$R$\textit{\ be a semiring with 
}$1\notin Z(R)$\textit{, and }$I$\textit{\ a left ideal of }$R$\textit{.
Then, the following statements are equivalent:}

\textit{(i) }$I\subseteq J(R)$\textit{;}

\textit{(ii) For any }$w$\textit{-finitely generated left }$R$\textit{%
-semimodule }$M$\textit{, if }$IM=M$\textit{, then }$M=Z(M)$\textit{%
.\medskip }

\noindent \textbf{Proof.} Because $1\notin Z(R)$, we first note that $%
R^{\ast }$and $D(R^{\ast })$ are nonzero additively cancellative semiring
and ring, respectively, with the identity$1^{\ast }$.

(i) $\Longrightarrow $ (ii). Let $IM=M$ for a left ideal $I\subseteq
J(R)\subseteq R$ of a semiring $R$ and an $w$-finitely generated left $R$%
-semimodule $M$, and $I^{\ast }:=\{r^{\ast }\,|\,r\in I\}\subseteq R^{\ast }$%
. It is clear that $I^{\ast }$ is a left ideal of $R^{\ast }$ and $I^{\ast
}\subseteq J(R^{\ast })$ since by \cite[p. 420, Section 4 e)]{i:otjroas} $%
J(R)=\{r\,|\,r^{\ast }\in J(R^{\ast })\}$. Also, it is clear that $D(I^{\ast
}):=\{r^{\ast }-s^{\ast }\,|\,r^{\ast },s^{\ast }\in I^{\ast }\}$ is a left
ideal of $D(R^{\ast })$, and even $D(I^{\ast })\subseteq J(D(R^{\ast }))$
because of $I^{\ast }\subseteq J(R^{\ast })$ and $J(D(R^{\ast }))$ is an
ideal of $D(R^{\ast })$, and, by \cite[p. 420, Section 4 e)]{i:otjroas}, $%
J(R^{\ast })=J(D(R^{\ast }))\cap R^{\ast }$. From $IM=M$, it follows that $%
I^{\ast }M^{\ast }=M^{\ast }$ and, hence, $D(I^{\ast })D(M^{\ast })$ $%
=D(M^{\ast })$. Applying Proposition 4.2 (iii) and (iv), one has that $%
D(M^{\ast })$ is a finitely generated left $D(R^{\ast })$-module. From these
observations and using the Nakayama's Lemma for rings (see, for example, %
\cite[Nakayama's Lemma 4.22]{lam:afcinr}), we have $D(M^{\ast })=0$; hence, $%
M^{\ast }=0$, and therefore, $M=Z(M)$.

(ii) $\Longrightarrow $ (i). In the notations introduced in the previous
implication and assuming that $I\nsubseteq J(R)$, we shall show that $%
D(I^{\ast })\nsubseteq J(D(R^{\ast }))$: Indeed, if $D(I^{\ast })\subseteq
J(D(R^{\ast }))$, then $I^{\ast }\subseteq J(D(R^{\ast }))$ and, by \cite[p.
420, Section 4 e)]{i:otjroas}, $I^{\ast }\subseteq R^{\ast }\cap J(D(R^{\ast
}))=J(R^{\ast })$; whence, using that by \cite[p. 420, Section 4 e)]%
{i:otjroas} $J(R)=\{r\,|\,r^{\ast }\in J(R^{\ast })\}$, we have that $%
I\subseteq J(R)$. Thus, there exists an irreducible left $D(R^{\ast })$%
-module $M$ and $D(I^{\ast })M\neq 0$. As $M$ is irreducible, $D(I^{\ast
})M=M$, and since $(M,+)$ is a group, we also have $I^{\ast }M=D(I^{\ast
})M=M$. For $M$ is an irreducible left $D(R^{\ast })$-module, we have that $%
M $ is an irreducible left $R$-semimodule with respect to the multiplication
defined by $rm=r^{\ast }m$ for all $r\in R$ and $m\in M$. Then, applying
Proposition 4.2 (i), we have $M$ is a nonzero $w$-finitely generated left $R$%
-semimodule with $IM=I^{\ast }M=M$, but $M\neq 0=Z(M)$. Therefore, $%
I\subseteq J(R)$. \textit{\ \ \ \ \ \ }$_{\square }\medskip $

As a corollary of Theorem 4.3, we obtain the following semiring version of
the Hopkins's Theorem for rings (see, \textit{e.g.}, \cite[Lemma 4.5.8]%
{gw:rtor}, or \cite[Theorem 4.12]{lam:afcinr}) .\medskip

\noindent \textbf{Corollary 4.4.} \textit{Let} $R$ \textit{be a left
artinian semiring. Then, there exists a natural number} $n\in \mathbb{N}$ 
\textit{such that} $J(R)^{n}=Z(R)$.\medskip

\noindent \textbf{Proof.} First note that for any irreducible $R$-semimodule 
$M$ and $z\in Z(R)$, always $z\in (0:M)_{R}$: Indeed, if $z+x=x$ for some $%
x\in M$, then $zm+xm=xm$ for any $m\in M$, and hence, $zm=0$ and $z\in
(0:M)_{R}$. Whence, $Z(R)\subseteq J(R)$. Therefore, if $1\in Z(R)$, then it
is clear that $R=Z(R)$, and hence, $J(R)^{n}=Z(R)$ for any $n\in \mathbb{N}$.

Now let $1\notin Z(R)$. For $R$ is a left artinian semiring, $R^{\ast }$ is
also left artinian semiring; hence, in the filtration 
\begin{equation*}
J(R^{\ast })\supseteq J(R^{\ast })^{2}\supseteq J(R^{\ast })^{3}\supseteq ...%
\text{,}
\end{equation*}%
$J(R^{\ast })^{n}=J(R^{\ast })^{n+1}$ for some positive $n\in \mathbb{N}$.
If $J(R^{\ast })^{n}\neq 0$, then the set of all left ideals $I$ of $R^{\ast
}$ for which $J(R^{\ast })^{n}I\neq 0$ is not empty. And therefore, by %
\cite[Proposition 2.1]{knt:ossss}, there exists a minimal left $I$ ideal of $%
R^{\ast }$with $J(R^{\ast })^{n}I\neq 0$. Let $x\in I$ and $J(R^{\ast
})^{n}x\neq 0$. Then, for $J(R^{\ast })x\subseteq R^{\ast }x\subseteq I$ and 
$J(R^{\ast })^{n}(J(R^{\ast })x)=J(R^{\ast })^{n+1}x=J(R^{\ast })^{n}x\neq 0$%
, one has $J(R^{\ast })x=I=R^{\ast }x$, and hence, by Theorem 4.3, $%
J(R^{\ast })(R^{\ast }x)=R^{\ast }x=I=0$. Thus, $J(R^{\ast })^{n}=0$, and
therefore, $J(R)^{n}=Z(R)$ since by \cite[p. 420, Section 4 e)]{i:otjroas} $%
J(R)=\{r\,|\,r^{\ast }\in J(R^{\ast })\}$. \textit{\ \ \ \ \ \ }$_{\square
}\medskip $

Recall (see, \textit{e.g.}, \cite[Chapter 4]{lam:afcinr}), a subring $%
R\subseteq End(_{D}M)$ of the endomorphism ring $End(_{D}M)$ of a left
vector space $_{D}M$ over a division ring $D$ is said to be \emph{dense} if,
for any linearly independent elements $m_{1},...,m_{n}\in $ $M$ and any
elements $m_{1}^{\prime },...,m_{n}^{\prime }\in $ $M$, there exists $r\in R$
such that $m_{i}r=m_{i}^{\prime }$ for $i=1,2,\ldots ,n$. Our next result is
a hemiring analog of the Jacobson-Chevalley Density Theorem for primitive
rings (see, for instance, \cite[Theorem 4.5.3]{gw:rtor}).\medskip

\noindent \textbf{Theorem 4.5.} \textit{A hemiring} $R$ \textit{is primitive
iff} \textit{it} \textit{is semiisomorphic to an additively cancellative
hemiring} $S$ \textit{whose ring of differences} $D(S)$ \textit{is, in turn,
isomorphic to a dense subring of linear transformations of a vector space
over a division ring.}\medskip

\noindent \textbf{Proof.} $\Longrightarrow $. Let $M$ be a faithful
irreducible left $R$-semimodule $M$, and $\varphi :R\longrightarrow End(M,+)$
the natural hemiring homomorphism such that $\varphi (r)(m)=rm$ for all $%
r\in R$ and $m\in M$ and, hence, having $ker(\varphi ):=\varphi
^{-1}(0)=(0:M)_{R}=0$. So, there exist the natural congruence $\equiv
_{\varphi }$on $R$ induced by $\varphi $, the natural surjection $\pi
:R\longrightarrow R/\equiv _{\varphi }$ with\ $ker(\pi ):=\pi ^{-1}(%
\overline{0})=(0:M)_{R}=0$, \textit{i.e.}, $\pi $ is a semiisomorphism, and
the natural injection $\psi :R/\equiv _{\varphi }\longrightarrow End(M,+)$
such that $\psi (\overline{r})=\varphi (r)$ for all $r\in R$. For $(M,+)$ is
additively cancellative, both $End(M,+)$ and $R/\equiv _{\varphi }$ are
additively cancellative hemirings, too. If $S$ stays for the hemiring $%
R/\equiv _{\varphi }$, $M$ becomes an irreducible left $S$-semimodule $_{S}M$
with $\overline{r}m=rm$ for all $m\in M$ and $r\in R$. Then, $D(M)$ is an
irreducible left $D(S)$-module with $(0:D(M))_{D(S)}=\{\overline{r}-%
\overline{s}\,|\,(\overline{r}-\overline{s})m=0$ for all $m\in M\}$ $=\{%
\overline{r}-\overline{s}\,|\,rm=sm\ $for all $m\in M\}=\{\overline{r}-%
\overline{s}$ $|$ $r\equiv _{\varphi }s\}=0$. Whence it follows that $D(S)$
is a primitive ring. Applying \cite[Theorem 4.5.3]{gw:rtor}, we immediately
obtain that $D(S)$ is isomorphic to a dense subring of linear
transformations of a space over a division ring.

$\Longleftarrow $. Let $\pi :R\longrightarrow S$ be a semiisomorphism to an
additively cancellative hemiring $S$ with $D(S)$ to be isomorphic to a dense
subring of linear transformations of a vector space over a division ring. By %
\cite[Theorem 4.5.3]{gw:rtor}, $D(S)$ is a primitive ring, and therefore,
there exists a faithful irreducible left $D(S)$-module $M$ which is also a
faithful irreducible left $S$-semimodule. Then $M$ is obviously a left $R$%
-semimodule with $rm=\pi (r)m$ for all $r\in R$ and $m\in M$; and $%
(0:M)_{R}=\{r\in R$ $|$ $\pi (r)m=0\ $for all $m\in M\}=\{r$ $|$ $\pi (r)\in
(0:M)_{S}=0\}$ $=ker(\pi )=0$. For $_{S}M$ is irreducible, $_{R}M$ is also
irreducible and, therefore, $R$ is a primitive hemiring. \textit{\ \ \ \ \ \ 
}$_{\square }\medskip $

Combining Theorem 4.5 and Corollary 3.8, one obtains\medskip

\noindent \textbf{Corollary 4.6.} \textit{A hemiring} $R$ \textit{is }$J$%
\textit{-semisimple iff} \textit{it} \textit{is semiisomorphic to a
subdirect product of some additively cancellative hemirings} $S$ \textit{%
whose rings of differences} $D(S)$ \textit{are isomorphic to dense subrings
of linear transformations of vector spaces over division rings.}\medskip

A \emph{right congruence} on a hemiring $R$ is just a congruence on the
(right) regular semimodule $R_{R}$. And we say that a hemiring $R$ is
(right) \emph{congruence-artinian} iff the descending chain condition on
right congruences on $R$ is held. Our next observation, a corollary of
Theorem 4.5, is a hemiring version of the classical structure theorem for
artinian primitive ring (see, \textit{e.g.}, \cite[Proposition 4.5.4]%
{gw:rtor}, or \cite[Theorem 4.11.19]{lam:afcinr}): \medskip

\noindent \textbf{Corollary 4.7.} \textit{A primitive congruence-artinian
hemiring} $R$ \textit{is semiisomorphic to an additively cancellative
hemiring} $S$ \textit{whose ring of differences }$D(S)$\textit{, in turn,} 
\textit{is isomorphic to a full ring of linear transformations on a finite
dimensional vector space over a division ring.}\medskip

\noindent \textbf{Proof.} By Theorem 4.5, there exists a semiisomorphism $%
f:R\longrightarrow S$ to an additively cancellative hemirings $S$ such that $%
D(S)$ is isomorphic to a dense subring of linear transformations $End(_{F}V)$
of a vector space $V$ over a division ring $F$. We shall show that $%
dim(_{F}V)<\infty $. Indeed, suppose that $\{e_{1},e_{2},...,e_{n},...\}$ is
a basis of $_{F}V$ and consider the subspaces $V_{n}$ generated by $%
e_{1},...,e_{n}$ for $n=1,2,...$. Obviously, all $(0:V_{n})_{D(S)}$ are
right ideals of $D(S)$, and 
\begin{equation*}
(0:V_{1})_{D(S)}\supset (0:V_{2})_{D(S)}\supset ...\supset
(0:V_{n})_{D(S)}\supset ...
\end{equation*}%
is a strictly descending chain according to \cite[Lemma 4.5.2]{gw:rtor}. It
is clear that $(0:V_{n})_{D(S)}=\{s-s^{\prime }\in D(S)$ $|$ $s,s^{\prime
}\in S$ $\&$ $\forall x\in V_{n}(xs=xs^{\prime })\}$ and $\alpha
_{n}:=\{(s,s^{\prime })\ |$ $s,s^{\prime }\in S$ $\&$ $\forall x\in
V_{n}(xs=xs^{\prime })\}$, $n=1,2,...$, are a right congruences on $S$.
Therefore, there is a strictly descending chain of right congruences $\alpha
_{n}$:\ 
\begin{equation*}
\alpha _{1}\supset \alpha _{2}\supset ...\supset \alpha _{n}\supset ...\quad
\quad \quad
\end{equation*}%
So, there is the strictly descending chain of right congruences $\beta
_{n}:=\{(r,r^{\prime })$ $|$ $(f(r),$ $f(r^{\prime }))\in \alpha _{n}\}$, $%
n=1,2,...$, on a congruence-artinian hemiring $R$. Thus, $dim(_{F}V)<\infty $%
, and therefore, $D(S)=End(_{F}V)$. \textit{\ \ \ \ \ \ }$_{\square
}\medskip $

As was mentioned earlier and in connection with Problem 1, we conclude this
section with the following observation:\medskip

\noindent \textbf{Proposition 4.8.} $J_{s}(R)\subseteq J(R)$ \textit{for any
commutative or additively regular, in particular additively idempotent,
hemiring }$R$\textit{.\medskip }

\noindent \textbf{Proof.} Let $R$ be a commutative hemiring. By %
\cite[Theorem 4.9]{m:otrtfs}, $J(R)=\cap \{I\in \mathcal{SI}(R)$ $|$ $%
J(R/I)=0\}$, where every commutative $J$-semisimple hemiring $R/I$ is, by
Corollary 3.8, semiisomorphic to a subdirect product of commutative
primitive hemirings. Taking into consideration Theorem 4.5, we note that a
commutative hemiring $S$ is primitive iff it is semiisomorphic to an
additively cancellative commutative hemiring $H$ whose ring of differences $%
D(H)$ is a field. If the hemiring $H$ is itself a ring, then it is a field
and $J_{s}(H)=J(H)=0$. Otherwise, we shall show that $H$ is semiisomorphic
to the semifield $\mathbf{B}$.

First, $H$ is zerosumfree: Indeed, if $a+b=0$ for some $a$, $b$ $\in $ $H$
and $a\neq 0$, then there exist $x,y\in H$ such that in the field of
differences $D(H)$ we have $a(x-y)=1$, that is, $ax=1+ay$ and, hence, $%
0=ax+bx=1+(ay+bx)$, what implies a contradiction with that $H$ is not a
ring. Secondly, $H$ is entire since it is a cancellative hemiring and
therefore is a subhemiring of the field $D(H)$. Thus, we have the
semiisomorphism $f:H\longrightarrow \mathbf{B}$ given by $h\longmapsto 1$
for all nonzero $h\in H$. As was mentioned in Example 3.7, $J_{s}(\mathbf{B}%
)=0$, and therefore, by \cite[Theorem 3.7 (b)]{hebwei:scos}, $J_{s}(H)=0$,
too.

From these observations and using \cite[Theorem 3.7 (b)]{hebwei:scos} one
more time, we conclude that $J_{s}(R/I)=0$ for every commutative $J$%
-semisimple hemiring $R/I$. Then, applying \cite[Theorem 4.9]{m:otrtfs}
again, we obtain the needed inclusion:%
\begin{eqnarray*}
J(R) &=&\cap \{I\in \mathcal{SI}(R)|\text{ }J(R/I)=0\}\supseteq \\
\cap \{I &\in &\mathcal{SI}(R)|\text{ }J_{s}(R/I)=0\}=J_{s}(R)\text{.}
\end{eqnarray*}%
\noindent \noindent

The statement for additively regular hemirings can be established in the
similar way by using Corollary 3.9 and the obvious fact that $J(R)=J_{s}(R)$
for rings $R$. \textit{\ \ \ \ \ \ }$_{\square }\medskip $

\section{Radicals and Morita Equivalence of Semirings}

Following \cite{kat:thcos} and \cite[theorem 4.12]{kn:meahcos}, two
semirings $R$ and $S$ are said to be \textit{Morita equivalent} if the
semimodule categories $_{R}\mathcal{M}$ and $_{S}\mathcal{M}$ are
equivalent, \textit{i.e.}, there exist two additive functors $F:$ $_{R}%
\mathcal{M}\longrightarrow $$_{S}\mathcal{M}$ and $G:$ $_{S}\mathcal{M}%
\longrightarrow $$_{R}\mathcal{M}$ and natural isomorphisms $\eta
:GF\longrightarrow Id_{_{R}\mathcal{M}}$ and $\xi :FG\longrightarrow Id_{_{S}%
\mathcal{M}}$. The semirings $R$ and $S$ are Morita equivalent iff the
semimodule categories $\mathcal{M}_{R}$ and $\mathcal{M}_{S}$ are also
equivalent \cite[Theorem 4. 12]{kn:meahcos}. Recall \cite{kn:meahcos}, a
left semimodule $_{R}P\in |_{R}\mathcal{M}|$ is a \textit{generator} for the
category of left semimodules $_{R}\mathcal{M}$ if the regular semimodule $%
_{R}R\in |_{R}\mathcal{M}|$ is a retract of a finite direct sum $\oplus
_{i}P $ of the semimodule $_{R}P$; and a finitely generated projective
semimodule $_{R}P\in |_{R}\mathcal{M}|$ is called a \textit{progenerator}
for the category $_{R}\mathcal{M}$ if it is generator for $_{R}\mathcal{M}$.
Also, a left semimodule $_{R}P\in |_{R}\mathcal{M}|$ is a progenerator for
the category of left semimodules $_{R}\mathcal{M}$ iff its trace ideal $%
tr(P):=\sum_{f\text{ }\in \text{ }_{R}\mathcal{M}(_{R}P,_{R}R)}f(P)$
coincides with $R$, \textit{i.e.}, $tr(P)=R$ \cite[Proposition 3.9]%
{kn:meahcos}. Then, the semirings $R$ and $S$ are Morita equivalent iff
there exists a progenerator $_{R}P\in |_{R}\mathcal{M}|$ for $_{R}\mathcal{M}
$ such that the semiring $S$ and the endomorphism semiring $End(_{R}P)$ are
isomorphic semirings \cite[Definition 4.1 and Theorem 4.12]{kn:meahcos}. For
the reader's convenience, we also reproduce here Proposition 5.2 of \cite%
{knz:ososacs} that will prove to be useful in a sequence:\medskip\ 

\noindent \textbf{Proposition 5.1.} (\textit{cf. }\cite[Proposition 18.33]%
{lam:lomar}) \textit{For semirings} $R$ \textit{and} $S$\textit{,} \textit{%
the followings statements are equivalent: }

\textit{(i)} $R$ \textit{is Morita equivalent to} $S$\textit{;}

\textit{(ii)} $S\cong eM_{n}(R)e$ \textit{for some idempotent} $e$ \textit{%
in a matrix semiring} $M_{n}(R)$ \textit{such that} $%
M_{n}(R)eM_{n}(R)=M_{n}(R)$\textit{.} \textit{\medskip }

\noindent \textbf{Proof. }(i) $\Longrightarrow $ (ii). By \cite[Theorem 4.12]%
{kn:meahcos}, there exists a progenerator $_{R}P\in |_{R}\mathcal{M}|$ for $%
_{R}\mathcal{M}$ such that $S\cong End(_{R}P)$ as semirings. Applying %
\cite[Proposition 3.1]{kn:meahcos} and without loss of generality, we can
assume that the semimodule $_{R}P$ is a subsemimodule of a free semimodule $%
_{R}R^{n}$, and there exists an endomorphism $e\in End(_{R}R^{n})$ such that 
$e^{2}=e$, $P=e(R^{n})$ and $e|_{P}=id_{P}$. For $e\in End(_{R}R^{n})\cong
M_{n}(R)$, one can consider $e$ to be a right multiplication by some
idempotent matrix $(a_{ij})\in M_{n}(R)$. In the same fashion as it has been
done in the case of modules over rings (see, for example, \cite[Remark
18.10(D) and Exercise 2.18]{lam:lomar}), we obtain that $tr(P)=\sum Ra_{ij}R$
and $rE_{ij}eE_{kl}r^{\prime }=ra_{jk}r^{\prime }E_{il}$, where $\{E_{ij}\}$
are the matrix units in $M_{n}(R)$ and $r,r^{\prime }\in R$, and get that $%
M_{n}(R)eM_{n}(R)=M_{n}(tr(R))$. For $_{R}P$ is a progenerator of the
category of semimodules $_{R}\mathcal{M}$ and \cite[Proposition 3.9]%
{kn:meahcos}, $tr(P)=R$ and, hence, $M_{n}(R)eM_{n}(R)=M_{n}(R)$. The proof
of the implication is now completed by noting that the semiring homomorphism 
\begin{equation*}
\theta :End(_{R}P)\longrightarrow eEnd(_{R}R^{n})e\cong eM_{n}(R)e
\end{equation*}%
defined for all $f\in End(_{R}P$ by $f\longmapsto eife$, where $%
i:P\longrightarrow R^{n}$ is the natural injection, is a semiring
isomorphism.

(ii) $\Longrightarrow $ (i). Let $S\cong eM_{n}(R)e$ for some idempotent $e$
of a matrix semiring $M_{n}(R)$ and $M_{n}(R)eM_{n}(R)=M_{n}(R)$. Then,
modifying for semimodules the well-known for modules over rings results
(see, for example, \cite[Proposition 21.6 and Corollary 21.7]{lam:afcinr}),
one has $S\cong eM_{n}(R)e\cong $

\noindent $End(_{M_{n}(R)}M_{n}(R)e)$ with $M_{n}(R)e$ to be, by %
\cite[Corollary 3.3]{kn:meahcos}, a projective left $M_{n}(R)$-semimodule.
Moreover, for $tr(_{M_{n}(R)}M_{n}(R)e)=M_{n}(R)eM_{n}(R)$

\noindent $=M_{n}(R)$, the semimodule $_{M_{n}(R)}M_{n}(R)e$ is a
progenerator of the category of semimodules $_{M_{n}(R)}\mathcal{M}$. From
these observations, it follows that $S$ and $M_{n}(R)$ are Morita equivalent
semirings, and therefore, applying \cite[Theorem 5.14]{kat:thcos} and %
\cite[Corollary 4.4]{kn:meahcos}, we conclude that the semirings $R$ and $S$
are Morita equivalent too. \textit{\ \ \ \ \ \ }$_{\square }\medskip $

A semiring $R$ is \textit{zeroic} \cite[p. 50]{golan:sata} if $Z(R)=R$, 
\textit{i.e.}, $1+r=r$ for some $r\in R$. From Proposition 5.1, it
immediately follows that to be zeroic is a Morita property for semirings,
namely:\medskip

\noindent \textbf{Corollary 5.2.} \textit{Let a zeroic semiring }$R$\textit{%
\ be Morita equivalent to a semiring }$S$\textit{. Then, }$S$\textit{\ is
also a zeroic semiring.\medskip }

\noindent \textbf{Proof. }It is clear that a matrix semiring $M_{n}(R)$ is
zeroic and, hence, the statement follows from Proposition 5.1 right away. 
\textit{\ \ \ \ \ \ }$_{\square }\medskip $

Let $D(\mathcal{M}_{R})$ denote the full subcategory of $R$\textit{-modules}
of the semimodule category $\mathcal{M}_{R}$ over a semiring $R$ with an $R$%
-semimodule $M\in |D(\mathcal{M}_{R})|$ iff the monoid $(M,+,0)$ is an
abelian group. For a nonzeroic semiring $R$, it is easy to see (see also %
\cite[Section 16, p. 183]{golan:sata}) that the category $D(\mathcal{M}_{R})$
actually coincides with the category $\mathcal{M}_{D(R)}$ of right modules
over the ring of differences $D(R)$. \medskip

\noindent \textbf{Proposition 5.3.} \textit{Let} $F:$ $\mathcal{M}%
_{R}\rightleftarrows \ \mathcal{M}_{S}:G$ \textit{be an equivalence of the
semimodule categories }$\mathcal{M}_{R}$ \textit{and }$\mathcal{M}_{S}$ 
\textit{over nonzeroic semirings }$R$\textit{\ and }$S$\textit{,
respectively. Then, for the same functors }$F$\textit{\ and }$G$\textit{,
the following statements are true:}

\textit{(i)} $F:$ $D(\mathcal{M}_{R})\rightleftarrows \ D(\mathcal{M}_{S}):G$
\textit{is an equivalence of the modules categories }$D(\mathcal{M}_{R})$ 
\textit{and} $D(\mathcal{M}_{S})$ \textit{over the semirings }$R$\textit{\
and }$S$\textit{, respectively;}

\textit{(ii)} $F:$ $\mathcal{M}_{D(R)}\rightleftarrows \ \mathcal{M}%
_{D(S)}:G $ \textit{is an equivalence of the modules categories }$\mathcal{M}%
_{D(R)}$ \textit{and }$\mathcal{M}_{D(S)}$ \textit{over the rings of
differences }$D(R)$\textit{\ and }$D(S)$\textit{, respectively;}

\textit{(iii) For any semimodules} $M\in |\mathcal{M}_{R}|$ \textit{and} $%
N\in |\mathcal{M}_{S}|$\textit{, there are semimodule isomorphisms of the
modules of differences }$F(D(M))_{S}\cong D(F(M))_{S}$\textit{\ and }$%
G(D(N))_{R}\cong D(G(N))_{R}$\textit{; in other words, the functors }$F$%
\textit{\ and }$G$\textit{\ preserve modules of differences;}

\textit{(iv)} \textit{The functors }$F$\textit{\ and }$G$\textit{\ preserve
the classes of cancellative semimodules.\medskip }

\noindent \textbf{Proof. }(i). This statement follows straightforwardly
noting that by \cite[Theorems 4.5 and 4.12]{kn:meahcos} the functor $F$, for
example, is isomorphic to additive functors $-\otimes _{R}P$ for a suitable
bisemimodule $_{R}P_{S}$.

(ii). It follows from item (i) and the remark above.

(iii). Let $M\in |\mathcal{M}_{R}|$. Then, using \cite[Proposition 16.1]%
{golan:sata}, it is easy to see that the $R$-module of differences $D(M)$ is
a colimit of some small diagram of semimodule surjections $\{$ $%
M\twoheadrightarrow X$ $|$ $X\in |D(\mathcal{M}_{R})|$ $\}$. By %
\cite[Theorem 5.5.1 and its dual]{macl:cwm}, the functor $F$ preserves
colimits and, therefore, now our assertions follows from (i).

(iv). As was mentioned in \cite[p. 181]{golan:sata}, a semimodule $M$ is
cancellative iff the canonical homomorphism $\xi _{M}:$ $M$ $\longrightarrow
D(M)$ given by $m\longmapsto \overline{(m,0)}$ (see Section 4) is injective.
From this observation, dual of \cite[Lemma 4.7]{kn:meahcos} and (iii), we
immediately deduce the statements. \textit{\ \ \ \ \ }$_{\square }\medskip $

We reformulate item (ii) of Proposition 5.3. as \medskip

\noindent \textbf{Corollary 5.4.} \textit{The Morita equivalence of
nonzeroic semirings }$R$\textit{\ and }$S$\textit{\ implies the Morita
equivalence of their rings of differences }$D(R)$\textit{\ and }$D(S)$%
\textit{. \medskip }

However, as the following obvious counterexample shows the inverse of this
statement, in general, is not true.\medskip

\noindent \textbf{Example 5.5. } \textit{The ring of integers} $\mathbb{Z}$ 
\textit{and the semiring} $\mathbb{N}$ \textit{are clearly not Morita
equivalent semirings, but} $D(\mathbb{N})=D(\mathbb{Z})=\mathbb{Z}$\textit{%
.\medskip }

Our next observation is that the classes of simple, semi-irreducible, and
irreducible, semimodules are also ``preserved'' by Morita equivalences of
semirings, more precisely:\medskip

\noindent \textbf{Proposition 5.6.} \textit{Let} $F:$ $\mathcal{M}%
_{R}\rightleftarrows \ \mathcal{M}_{S}:G$ \textit{be an equivalence of the
semimodule categories }$\mathcal{M}_{R}$ \textit{and }$\mathcal{M}_{S}$ 
\textit{over semirings }$R$\textit{\ and }$S$\textit{, respectively. Then, a
semimodule }$M\in |\mathcal{M}_{R}|$ \textit{is simple (semi-irreducible,
irreducible) iff }$F(M)\in |\mathcal{M}_{S}|$ \textit{is simple
(semi-irreducible, irreducible).\medskip }

\noindent \textbf{Proof. }First, note the following obvious facts for a
nonzero semimodule $M\in |\mathcal{M}_{R}|$: a) $M$ has no nonzero
subsemimodules iff every nonzero homomorphism $f:L\longrightarrow M$ in $%
\mathcal{M}_{R}$ is a surjection; b) $M$ is congruence-simple iff every
nonzero homomorphism $g:M\longrightarrow N$ in $\mathcal{M}_{R}$ is an
injection; c) a cancellative semimodule $M$ is semi-irreducible iff $%
Ker(f)=0 $ for every nonzero homomorphism $f:M\longrightarrow N$ in $%
\mathcal{M}_{R}$.

Now, let $M\in |\mathcal{M}_{R}|$ be a simple semimodule, and $%
f:L\longrightarrow F(M)$ a nonzero homomorphism in $\mathcal{M}_{S}$. Then, $%
G(f):G(L)\longrightarrow G(F(M))\cong M$ is also a nonzero homomorphism in $%
\mathcal{M}_{R}$, and therefore, it is a surjection. Whence by \cite[Lemma
4.7]{kn:meahcos}, $F(G(f)):F(G(L))\longrightarrow F(G(F(M)))\cong F(M)$ and,
hence, $f$ are also surjections and $F(M)$ has no nonzero subsemimodules.
For a nonzero homomorphism $g:F(M)\longrightarrow N$ in $\mathcal{M}_{S}$,
we have $G(g):M\cong G(F(M))\longrightarrow G(N)$ is a nonzero homomorphism
in $\mathcal{M}_{R}$, and hence, it is injective. Then, applying the functor 
$F$ and dual of \cite[Lemma 4.7]{kn:meahcos}, we have that $%
F(G(g)):F(M)\cong F(G(F(M)))\longrightarrow F(G(N))\cong N$ is an injection
and, hence, $g$ is also an injection. Thus, we have shown that $F(M)$ is a
simple semimodule too.

Let $M\in |\mathcal{M}_{R}|$ be a semi-irreducible semimodule, and $%
g:F(M)\longrightarrow N$ a nonzero homomorphism in $\mathcal{M}_{S}$. Then, $%
G(g):M\cong G(F(M))\longrightarrow G(N)$ is a nonzero homomorphism in $%
\mathcal{M}_{R}$ with $Ker(G(g))=0$. Again, it is easy to see that $%
Ker(G(g)) $ is a limit of some small diagram of semimodule injections $%
\{i:X\rightarrowtail G(F(M))$ $|$ $X\in |\mathcal{M}_{R}|$ $\&$ $G(g)i=0\}$,
and therefore, applying \cite[Theorem 5.5.1 ]{macl:cwm}, we have $%
Ker(F(G(g)):F(M)\cong F(G(F(M)))\longrightarrow F(G(N)))\cong N)=0$, and
hence, $Ker(g)=0$ and $F(M)\in |\mathcal{M}_{S}|$ is a semi-irreducible
semimodule.

Finally, consider the case when a semimodule $M\in |\mathcal{M}_{R}|$ is
irreducible. By Corollary 5.2, $R$ and $S$ are zeroic semirings
simultaneously; and, clearly, over zeroic semirings there are no irreducible
semimodules. Hence, for zeroic semirings $R$ and $S$ our assertion is true.
So, assume that the semirings $R$ and $S$ are not zeroic. As was mentioned
earlier, $M\in |\mathcal{M}_{R}|$ is irreducible iff $D(M)\in |\mathcal{M}%
_{D(R)}|$ is irreducible, and, applying Proposition 5.3 , iff $F(D(M))\in |%
\mathcal{M}_{D(S)}|$ is irreducible, iff $D(F(M))\in |\mathcal{M}_{D(S)}|$
is irreducible and, finally, iff $F(M)\in |\mathcal{M}_{S}|$ is irreducible. 
\textit{\ \ \ \ \ \ }$_{\square }\medskip $

In our next result, we establish the fundamental relationship between the
radicals $J$, $J_{s}$, and $\mathcal{R}_{BM}$ of semirings $R$ and matrix
semirings $M_{n}(R)$, $n\geq 1$. But first note the following fact: \medskip

\noindent \textbf{Lemma 5.7.} \textit{A semiring} $R$ \textit{is }$\mathcal{R%
}_{BM}$\textit{-semisimple iff all matrix semirings }$M_{n}(R)$\textit{,} $%
n\geq 1$\textit{, are }$\mathcal{R}_{BM}$\textit{-semisimple, i.e., }$%
\mathcal{R}_{BM}(R)=0$ \textit{iff} \textit{\ }$\mathcal{R}_{BM}(M_{n}(R))=0$
\textit{for all} $n\geq 1$\textit{.\medskip }

\noindent \textbf{Proof.} In the same way as for matrix rings over rings, it
is easy to show that any ideal of a matrix semiring $M_{n}(R)$ over a
semiring $R$ consists of all matrices having all entries in some ideal of $R$%
. Hence, for a nonzero ideal $A$ of the semiring $M_{n}(R)$ there exists a
nonzero ideal $I$ of a semiring $R$ such that $A=M_{n}(I)$. If $\ $the
semiring $R$ is $\mathcal{R}_{BM}$-semisimple, \textit{i.e.}, $\mathcal{R}%
_{BM}(R)=0$, then $I\nsubseteq \mathcal{R}_{BM}(R)$ and there exists a
nonzero surjective semiring homomorphism $f$ from $I$ to an ideal-simple
semiring $S$. Whence, there exists a nonzero surjective homomorphism from $%
A=M_{n}(I)$ to the matrix semiring $M_{n}(S)$ that, by \cite[Proposition 4.7]%
{knt:mosssparp}, is also an ideal-simple semiring. Therefore, $\mathcal{R}%
_{BM}(M_{n}(R))=0$. \textit{\ \ \ \ \ \ }$_{\square }\medskip $

\noindent \textbf{Theorem 5.8.} \textit{For all matrix semirings }$M_{n}(R)$%
\textit{,} $n\geq 1$\textit{, over a semiring} $R$\textit{, the following
equations are held: }

\textit{(i)} $\mathcal{R}_{BM}(M_{n}(R))=M_{n}(\mathcal{R}_{BM}(R))$\textit{;%
}

\textit{(ii)} $J_{s}(M_{n}(R))=M_{n}(J_{s}(R))$\textit{;}

\textit{(iii)} $J(M_{n}(R))=M_{n}(J(R))$\textit{.\medskip }

\noindent \textbf{Proof.} (i). For $\mathcal{R}_{BM}(R)$ is a subtractive
ideal of $R$, it follows that $M_{n}(\mathcal{R}_{BM}(R))$ is a subtractive
ideal of $M_{n}(R)$ and, hence,

\noindent $M_{n}(R)/M_{n}(\mathcal{R}_{BM}(R))\cong M_{n}(R/\mathcal{R}%
_{BM}(R))$. So, because of Lemma 5.2 and $\mathcal{R}_{BM}(R/\mathcal{R}%
_{BM}(R))=0$, we get $\mathcal{R}_{BM}(M_{n}(R)/M_{n}(\mathcal{R}_{BM}(R)))=$

\noindent $\mathcal{R}_{BM}(M_{n}(R/\mathcal{R}_{BM}(R)))=0$. From the
latter and since by \cite[Theorem 4.9]{m:otrtfs} $\mathcal{R}%
_{BM}(M_{n}(R))=\cap \{K\in \mathcal{SI}(M_{n}(R))\,|\,\mathcal{R}%
_{BM}(M_{n}(R)/K)=0\}\subseteq M_{n}(\mathcal{R}_{BM}(R))$, we have an
inclusion $\mathcal{R}_{BM}(M_{n}(R))\subseteq $ $M_{n}(\mathcal{R}_{BM}(R))$%
.

Let $I:=\mathcal{R}_{BM}(R)$, and $\mathcal{R}_{BM}(M_{n}(I))\neq M_{n}(I)$.
So, there exists a congruence $\rho $ on $M_{n}(I)$ such that $M_{n}(I)/\rho 
$ is an ideal-simple semiring with an identity $\overline{E}$ for $E\in
M_{n}(I)$. By Lemma 3.14, $\overline{\rho }:=\{(A,B)\in M_{n}(R)^{2}\,|\,$

\noindent $(EAE,EBE)\in \rho \}$ is a congruence on $M_{n}(R)$, and $\varphi
:M_{n}(R)/\overline{\rho }\longrightarrow M_{n}(I)/\rho $, given by $%
\overline{A}\longmapsto \overline{EAE}$, is a semiring isomorphism. It is
easy to see that the relation $\theta $ defined for $a,b\in R$ as follows: 
\begin{equation*}
a\theta b\Longleftrightarrow \forall i,j\text{ }(aE_{ij}\overline{\rho }%
bE_{ij})\text{, where }\{E_{ij}\}\text{ are the matrix units in }M_{n}(R)%
\text{,}
\end{equation*}%
is a congruence on $R$ and $M_{n}(R/\theta )\cong M_{n}(R)/\overline{\rho }$%
. Then for $\delta :=\theta \cap J^{2}$, one can readily see that $\delta $
is a congruence on $I$ and there is the injective homomorphism $\psi
:I/\delta \longrightarrow R/\theta $ given by $I/\delta \ni \overline{a}%
\longmapsto \overline{a}\in R/\theta $ for all $a\in I$, and $Im(\psi )$ is
an ideal of $R/\theta $. For $M_{n}(I)/\rho \neq 0$, there exists an element 
$a\in I$ such that $(aE_{ij},0)\notin \rho $ for some $i,j$, and hence, $%
(aE_{ij},0)\notin \overline{\rho }$. Therefore, $(a,0)\notin \delta $, and $%
Im(\psi )$ is a nonzero ideal of $R/\theta $. Since $M_{n}(I)/\rho $ is
ideal-simple, so is $M_{n}(R/\theta )$, and by \cite[Proposition 4.7]%
{knt:mosssparp}, $R/\theta $ is also an ideal-simple semiring. Whence, $%
Im(\psi )=R/\theta $, that is, $\psi $ is an isomorphism, and hence, $%
I/\delta $ is an ideal-simple semiring in a contradiction with $I=\mathcal{R}%
_{BM}(R)$. Thus, $\mathcal{R}_{BM}(M_{n}(I))=M_{n}(I)$, and hence, $%
M_{n}(I)\subseteq \mathcal{R}_{BM}(M_{n}(R))$. Therefore, $M_{n}(\mathcal{R}%
_{BM}(R))=\mathcal{R}_{BM}(M_{n}(R))$.

(ii). First note that a simple semimodule $M\in |\mathcal{M}_{R}|$ is always
unitary, that is, $1.m=m$ for all $m\in M$. Indeed, $r(1.m)=(r.1)m=rm$ for
all $r\in R$ and consider the congruence $\beta $ on $M$ given by: $(x,y)\in
\beta $ iff $rx=ry$ for all $r\in R$ and $x,y\in M$. Since $M$ is
congruence-simple, $\beta =id_{M}$ or $\beta =M^{2}$. If $\beta =M^{2}$,
then $(x,0)\in \beta $, for all $x\in M$, that is, $rx=0$ for all $r\in R$
and $x\in M$, and hence, $RM=0$. Thus, $\beta =id_{M}$, that is, $(1.m,m)\in
\beta $,\textit{\ i.e.}, $1.m=m$.

By \cite[Theorem 5.14]{kat:thcos} and \cite[Theorem 4.12]{kn:meahcos}, the
semirings $R$ and $M_{n}(R)$ are Morita equivalent via the equivalence $F:$ $%
\mathcal{M}_{R}\rightleftarrows \ \mathcal{M}_{M_{n}(R)}:G$ such that $%
F(A)=A^{n}$ and $G(B)=E_{11}B$ for all $A\in |\mathcal{M}_{R}|$ and $B\in |%
\mathcal{M}_{M_{n}(R)}|$, where $E_{11}$ denotes the matrix unit. By
Proposition 5.6, $M^{n}=F(M)$ is also a simple $M_{n}(R)$-semimodule, and
therefore, denoting by $\mathcal{S}_{R}$ the set of all simple $R$%
-semimodules, we obtain an inclusion 
\begin{equation*}
J_{s}(M_{n}(R))\subseteq \cap _{M\in \mathcal{S}_{R}}(0:M^{n})_{M_{n}(R)}%
\subseteq M_{n}(J_{s}(R))\text{.}
\end{equation*}

Similarly, for a simple semimodule $A\in \mathcal{M}_{M_{n}(R)}$, by
Proposition 5.6, $E_{11}A=G(A)$ is also a simple left $R$-semimodule, and
therefore, denoting by $\mathcal{S}_{M_{n}(R)}$ the set of all simple $%
M_{n}(R)$-semimodules and noting the obvious inclusion $J_{s}(R)\subseteq
\cap _{A\in \mathcal{S}_{M_{n}(R)}}(0:E_{11}A)_{R}$, we have the opposite
inclusion 
\begin{equation*}
M_{n}(J_{s}(R))\subseteq M_{n}(\cap _{A\in \mathcal{S}%
_{M_{n}(R)}}(0:E_{11}A)_{R})\subseteq J_{s}(M_{n}(R)).
\end{equation*}

(iii). Just using the fact that $J(R)=\{(0:M)_{R}|$ $M\in $ $|\mathcal{M}%
_{R}|$ is irreducible$\}$, this equation can be proved similarly to (ii).
Also, it was established, but by a different method, in \cite[Theorem 9]%
{b:tjroas}. \textit{\ \ \ \ \ \ }$_{\square }\medskip $

Moreover, we shall show that Theorem 5.8 can be extended to hemirings in
general, but for that we first need some useful facts.\medskip

\noindent \textbf{Lemma 5.9.} \textit{Let} $\mathbb{R}$ \textit{be a radical
class of} $\mathbb{H}$\textit{,} \textit{and} $\varrho _{\mathbb{R}}$ 
\textit{its radical operator having} $\varrho _{\mathbb{R}}(\mathbb{N})=0$%
\textit{.} \textit{Then, }$\varrho _{\mathbb{R}}(R)=\varrho _{\mathbb{R}%
}(R^{1})$\textit{\ for any hemiring} $R$\textit{.\medskip }

\noindent \textbf{Proof.} It is clear that $R$ is a subtractive ideal of $%
R^{1}$ and $R^{1}/R\cong \mathbb{N}$. Therefore, $\varrho _{\mathbb{R}%
}(R^{1}/R)=\varrho _{\mathbb{R}}(\mathbb{N})=0$. By \cite[Theorem 4.9]%
{m:otrtfs}, $\varrho _{\mathbb{R}}(R^{1})=\cap \{K\in \mathcal{SI}%
(R^{1})\,|\,\varrho _{\mathbb{R}}(R^{1}/K)=0\}$, and hence, $\varrho _{%
\mathbb{R}}(R^{1})\subseteq R$. From this observation and \cite[Theorem 6.2]%
{m:otrtfs} (the \textit{A-D-S-property} of a radical class), we have $%
\varrho _{\mathbb{R}}(R^{1})\subseteq \varrho _{\mathbb{R}}(R)$. On the
other hand, because of $R$ is an ideal of $R^{1}$ and \cite[Theorem 6.2]%
{m:otrtfs}, $\varrho _{\mathbb{R}}(R)\subseteq \varrho _{\mathbb{R}}(R^{1})$%
, and hence, $\varrho _{\mathbb{R}}(R^{1})=\varrho _{\mathbb{R}}(R).$\textit{%
\ \ \ \ \ \ }$_{\square }\medskip $

\noindent \textbf{Lemma 5.10.} \textit{\ }$J(\mathbb{N})=\mathcal{R}_{BM}(%
\mathbb{N})=$ $J_{s}(\mathbb{N})=0$\textit{.\medskip }

\noindent \textbf{Proof.} For all prime numbers $p$, one can readily see
that $\mathbb{N}/p\mathbb{N}=\mathbb{Z}_{p}$ is a irreducible $\mathbb{N}$%
-semimodule. Therefore, $J(\mathbb{N})\subseteq \cap \{(0:\mathbb{Z}_{p})_{%
\mathbb{N}}\,|$ $p$ is prime$\}=\cap \{p\mathbb{N}\,|$ $p$ is prime$\}=0$,
and hence, $J(\mathbb{N})=0$.

(ii). $\mathcal{R}_{BM}(\mathbb{N}/p\mathbb{N})=\mathcal{R}_{BM}(\mathbb{Z}%
_{p})=0$ and $p\mathbb{N}$ is a subtractive ideal of $\mathbb{N}$ for all
prime numbers $p$, and hence, by \cite[Theorem 4.9]{m:otrtfs}, $\mathcal{R}%
_{BM}(\mathbb{N})=\cap \{K\in \mathcal{SI}(\mathbb{N})\,|\,\mathcal{R}_{BM}(%
\mathbb{N}/K)=0\}\subseteq \cap \{p\mathbb{N}\,|$ $p$ is prime$\}=0$.

(iii). Obviously, $(0:\mathbf{B})_{\mathbb{N}}=0$ for the simple $\mathbb{N}$%
-semimodule $\mathbf{B}_{\mathbb{N}}$ $\in $ $|\mathcal{M}_{\mathbb{N}}|$,
and therefore, $J_{s}(\mathbb{N})=0$. \textit{\ \ \ \ \ \ }$_{\square }$
\medskip

\noindent \textbf{Corollary 5.11.} \textit{For all matrix hemirings }$%
M_{n}(R)$\textit{,} $n\geq 1$\textit{, over a hemiring} $R$\textit{, the
following equations are held: }

\textit{(i)} $\mathcal{R}_{BM}(M_{n}(R))=M_{n}(\mathcal{R}_{BM}(R))$\textit{;%
}

\textit{(ii)} $J_{s}(M_{n}(R))=M_{n}(J_{s}(R))$\textit{;}

\textit{(iii)} $J(M_{n}(R))=M_{n}(J(R))$\textit{.}$\medskip $

\noindent \textbf{Proof.} (i). Applying Theorem 5.8 and Lemmas 5.9 and 5.10,
one gets $M_{n}(\mathcal{R}_{BM}(R))=M_{n}(\mathcal{R}_{BM}(R^{1}))=\mathcal{%
R}_{BM}(M_{n}(R^{1}))$. Also, since $\mathcal{R}_{BM}(R^{1}/R)$

\noindent $=\mathcal{R}_{BM}(\mathbb{N})=0$ and Lemma 5.7, one has $\mathcal{%
R}_{BM}(M_{n}(R^{1})/M_{n}(R))\cong $

\noindent $\mathcal{R}_{BM}(M_{n}(R^{1}/R))=0$. Then, noting that $M_{n}(R)$
is a subtractive ideal of $M_{n}(R^{1})$, since $R$ is a subtractive ideal
of $R^{1}$, and applying \cite[Theorem 4.9]{m:otrtfs} and \cite[Theorem 6.2]%
{m:otrtfs}, we have the inclusions $\mathcal{R}_{BM}(M_{n}(R^{1}))=\cap
\{K\in \mathcal{SI}(M_{n}(R^{1}))\,|\,\mathcal{R}_{BM}(M_{n}(R^{1})/K)=0\}%
\subseteq M_{n}(R)$ and $\mathcal{R}_{BM}(M_{n}(R^{1}))\subseteq \mathcal{R}%
_{BM}(M_{n}(R))$. Furthermore, noting that $M_{n}(R)$ is an ideal of $%
M_{n}(R^{1})$, since $R$ is an ideal of $R^{1}$, and applying \cite[Theorem
6.2]{m:otrtfs} again, we have the opposite inclusion $\mathcal{R}%
_{BM}(M_{n}(R^{1}))\supseteq \mathcal{R}_{BM}(M_{n}(R))$, and therefore, $%
\mathcal{R}_{BM}(M_{n}(R))$

\noindent $=M_{n}(\mathcal{R}_{BM}(R))$.

Two other equations, (ii) and (iii), can be justified in the similar
fashion. \textit{\ \ \ \ \ \ }$_{\square }\medskip $

As for rings (see, for example, \cite[Section 4.9]{gw:rtor}), we say that a
class $\delta $ of hemirings is \emph{matric-extensible}, if for all natural
numbers $n$, $\ $a hemiring $R\in \delta $ iff $M_{n}(R)\in \delta $.
\medskip

\noindent \textbf{Lemma 5.12.} \textit{For any} \textit{radical operator }$%
\varrho $ \textit{on} $\mathbb{H}$\textit{, hemiring} $R$ \textit{and
natural number} $n$\textit{,} \textit{there exists an ideal }$I$\textit{\ of 
}$R$\textit{\ such that} $\varrho (M_{n}(R))=M_{n}(I)$\textit{.}\medskip

\noindent \textbf{Proof.} For $M_{n}(R)$ is an ideal of the semiring $%
M_{n}(R^{1})$ and \cite[Theorem 6.2]{m:otrtfs}, $\varrho (M_{n}(R))$ is an
ideal of $M_{n}(R^{1})$. So there exists an ideal $I$ of $R^{1}$ such that $%
\varrho (M_{n}(R))=M_{n}(I)$; and since $M_{n}(I)\subseteq M_{n}(R)$, one
gets $I\subseteq R$. \textit{\ \ \ \ \ \ }$_{\square }\medskip $

Using this lemma, one readily obtains the following hemiring analog of %
\cite[Theorem 4.9.3]{gw:rtor}:\medskip

\noindent \textbf{Proposition 5.13.} \textit{For a radical class }$\gamma $ 
\textit{of hemirings, the following statements are equivalent:}

\textit{(i) }$\gamma $\textit{\ is matric-extensible;}

\textit{(ii) For every hemiring }$R\ $\textit{and natural number }$n$\textit{%
, a matrix equation }$\gamma (M_{n}(R))=M_{n}(\gamma (R))$\textit{\ is true;}

\textit{(iii) The semisimple class }$S_{\gamma }$\textit{\ of }$\gamma $%
\textit{\ is a matric-extensible class.\medskip }

\noindent \textbf{Proof.} (i) $\Longrightarrow $ (ii). Using Lemma 5.12 and
proceeding in a similar fashion as it has been done in the proof of the
implication (1) $\Longrightarrow $ (3) in \cite[Theorem 4.9.3]{gw:rtor}, one
gets\ this implication.

(ii) $\Longrightarrow $ (iii). This can be proved in a similar way as the
implication (3) $\Longrightarrow $ (2) in \cite[Theorem 4.9.3]{gw:rtor}.

(iii) $\Longrightarrow $ (i). Let $R\in \gamma $, but $M_{n}(R)\notin \gamma 
$. Then, $0\neq M_{n}(R)/\gamma (M_{n}(R))\in \mathcal{S}_{\gamma }$, and
from Lemma 5.12 it follows that $\gamma (M_{n}(R))=M_{n}(I)$ for some ideal $%
I$ of $R$. Hence, $0\neq M_{n}(R/I)\cong M_{n}(R)/M_{n}(I)=M_{n}(R)/\gamma
(M_{n}(R))\in \mathcal{S}_{\gamma }$. For $\mathcal{S}_{\gamma }$ is
matrix-extensible, $R/I\in \mathcal{S}_{\gamma }$; on the other hand, by
Theorem 2.2 (3), $R/I\in \gamma $. Thus,\ $R\in \gamma $ implies $%
M_{n}(R)\in \gamma $.

Now, let $M_{n}(R)\in \gamma $, but $R\notin \gamma $. Then, $\gamma
(R)\varsubsetneq R$ and because $\gamma (R)$ is a subtractive ideal of $R$
we have $0\neq R/\gamma (R)\in \mathcal{S}_{\gamma }$. Hence, $%
M_{n}(R)/M_{n}(\gamma (R))$

\noindent $\cong M_{n}(R/\gamma (R))\in \mathcal{S}_{\gamma }$ and $%
M_{n}(R/\gamma (R))\neq 0$; however, $M_{n}(R/\gamma (R))=0$ for $%
M_{n}(R)\in \gamma $. Thus, $M_{n}(R)\in \gamma $ implies $R\in \gamma $.%
\textit{\ \ \ \ \ \ }$_{\square }\medskip $

From Corollary 5.11 and Proposition 5.13, we immediately obtain \medskip

\noindent \textbf{Theorem 5.14.} \textit{The radical classes of hemirings of
the Jacobson radical }$J$\textit{, Brown-McCoy radical }$\mathcal{R}_{BM}$,%
\textit{\ and} \textit{the radical }$J_{s}$, \textit{are matrix-extensible
classes.\medskip }

Next we present some ``computational'' results regarding radicals of
hemirings $R$ and $eRe$ for idempotents $e\in R$, which, in our view, are
interesting and important on their own, and will prove to be useful in a
sequel.\medskip

\noindent \textbf{Proposition 5.15.} \textit{For any} \textit{idempotent }$e$
\textit{of a hemiring} $R$\textit{, the following statements are true:}

\textit{(i)} $J(eRe)=eJ(R)e$\textit{;}

\textit{(ii)} $eRe/J(eRe)\cong \overline{e}$$\overline{R}$$\overline{e}$%
\textit{,} \textit{where} $\overline{e}$ \textit{is the image of} $e$ 
\textit{in} $\overline{R}=R/J$\textit{.\medskip }

\noindent \textbf{Proof.} (i). For by \cite[Theorem 3]{b:tjroas} $J(R)$ is a
right semiregular ideal of $R$ (see, \cite[Definition 3]{b:tjroas}), $%
eae,ebe\in J(R)$ for any $a,b\in J(R)$ and there exist elements $%
r_{1},r_{2}\in J(R)$ such that 
\begin{equation*}
eae+r_{1}+eaer_{1}+eber_{2}=ebe+r_{2}+eaer_{2}+eber_{1}\text{.}
\end{equation*}%
So, as $e$ is an idempotent, 
\begin{equation*}
eae+er_{1}e+eae.er_{1}e+ebe.er_{2}e=ebe+er_{2}e+eae.er_{2}e+ebe.er_{1}e\text{%
.}
\end{equation*}%
Therefore, $eJ(R)e$ is a right semiregular ideal of $eRe$, and hence, by %
\cite[theorem 3 and Definition 4]{b:tjroas}, $eJ(R)e\subseteq J(eRe)$.

Now we shall show that $eJ(R)e\supseteq J(eRe)$. Indeed, one can readily see
that $(eRe)^{\ast }=\{r^{\ast }\in R^{\ast }\,|\,r\in eRe\}$, $D((eRe)^{\ast
})=\{r^{\ast }-s^{\ast }\,|\,r^{\ast },s^{\ast }\in (eRe)^{\ast }\}$, $%
(eRe)^{\ast }=e^{\ast }R^{\ast }e^{\ast }$, $D((eRe)^{\ast })=e^{\ast
}D(R^{\ast })e^{\ast }$ and, using \cite[Theorem 21.10]{lam:afcinr}, $%
J(e^{\ast }D(R^{\ast })e^{\ast })=e^{\ast }J(D(R^{\ast }))e^{\ast }$.
Therefore, 
\begin{equation*}
J((eRe)^{\ast })=(eRe)^{\ast }\cap J(D((eRe)^{\ast }))=
\end{equation*}%
\begin{equation*}
=(eRe)^{\ast }\cap J(e^{\ast }D(R^{\ast })e^{\ast })=(eRe)^{\ast }\cap
e^{\ast }J(D(R^{\ast }))e^{\ast }
\end{equation*}%
\begin{equation*}
\subseteq e^{\ast }(R^{\ast }\cap J(D(R^{\ast })))e^{\ast }=e^{\ast
}J(R^{\ast })e^{\ast }\subseteq J(R^{\ast }).
\end{equation*}

So, for any $x\in J(eRe)$, we have that $x^{\ast }\in J((eRe)^{\ast
})\subseteq J(R^{\ast })$, and hence, by \cite[Section 4 e), p. 420]%
{i:otjroas}, $x\in J(R)$; as $x\in J(eRe)\subseteq eRe$, we have $x=exe\in
eJ(R)e$, and therefore, $J(eRe)\subseteq eJ(R)e$.

(ii). Consider the natural hemiring homomorphism $f:eRe\longrightarrow $ $%
\overline{e}$$\overline{R}$$\overline{e}$ given by $ere\longmapsto \overline{%
e}\,\overline{r}\,\overline{e}$. It is clear that $f$ induces a surjection $%
g:eRe/J(eRe)$

\noindent $\twoheadrightarrow $ $\overline{e}$$\overline{R}$$\overline{e}$
that, in fact, is an isomorphism: Indeed, if $\overline{e}\,\overline{r}\,%
\overline{e}=$ $\overline{e}\,\overline{s}\,\overline{e}$, then there exist
elements $a,b\in J(R)$ such that $ere+a=ese+b$, and hence, $ere+eae=ese+ebe$
and, since by (i) $eae,ebe\in eJ(R)e=J(eRe)$, one has $ere\equiv
_{J(eRe)}ese $. \textit{\ \ \ \ \ \ }$_{\square }\medskip $

\noindent \textbf{Proposition 5.16.} \textit{\ (i) For any }$\mathcal{R}%
_{BM} $\textit{-semisimple} \textit{semiring} $R$ \textit{and an idempotent }%
$e\in $ $R$\textit{, the semiring }$eRe$ \textit{is also} $\mathcal{R}_{BM}$%
\textit{-semisimple.}

\textit{(ii) Let} $e\in R$ \textit{be a full idempotent, i.e. }$ReR=R$%
\textit{, in a hemiring} $R$\textit{.} \textit{Then,} $\mathcal{R}%
_{BM}(eRe)=e\mathcal{R}_{BM}(R)e$ \textit{and} $eRe/\mathcal{R}%
_{BM}(eRe)\cong \overline{e}$$\overline{R}$$\overline{e}$\textit{,} \textit{%
where} $\overline{e}$ \textit{is the image of} $e$ \textit{in} $\overline{R}%
=R/\mathcal{R}_{BM}(R)$\textit{.}

\textit{(iii)} \textit{Let} $e\in R$ \textit{be a full idempotent in a
semiring} $R$. \textit{Then,} $J_{s}(eRe)=eJ_{s}(R)e$ \textit{and} $%
eRe/J_{s}(eRe)\cong \overline{e}$$\overline{R}$$\overline{e}$\textit{,} 
\textit{where} $\overline{e}$ \textit{is the image of} $e$ \textit{in} $%
\overline{R}=R/J_{s}(R)$\textit{.\medskip }

\noindent \textbf{Proof.} (i). Let $I$ be a nonzero ideal of the semiring $%
eRe$, and $A:=RIR$ a nonzero ideal of $R$. For $\mathcal{R}_{BM}(R)=0$,
there exists a congruence $\rho $ on $A$ such that $A/\rho $ is an
ideal-simple semiring with identity $\overline{e}_{1}$. By Lemma 3.14, $\ 
\overline{\rho }:=\{(a,b)\in R^{2}\,|\,(e_{1}ae_{1},e_{1}be_{1})\in \rho \}$
is a congruence on $R$, and hence, $\theta =\overline{\rho }\cap I^{2}$ is a
congruence on $I$. We shall show that $I/\theta $ is an ideal-simple
semiring.

Indeed, for $e_{1}\in A$ there exist elements $%
r_{1},...,r_{k},s_{1},...,s_{k}\in R$ and $a_{1},...,a_{k}\in I$ such that $%
e_{1}=r_{1}a_{1}s_{1}+...+r_{k}a_{k}s_{k}$ with $ea_{i}e=a_{i}$ for all $%
a_{i}$, $i=1,\ldots ,k$; and let $%
e_{2}:=ee_{1}e=er_{1}ea_{1}es_{1}e+...+er_{k}ea_{k}es_{k}e\in I$. Then, $%
(ae_{2},a)\in I^{2}$, $ae=a$, and hence, $\overline{e_{1}ae_{2}e_{1}}=%
\overline{e_{1}aee_{1}ee_{1}}=\overline{e_{1}ae_{1}ee_{1}}=\overline{%
e_{1}aee_{1}}=\overline{e_{1}ae_{1}}$ and $(ae_{2},a)\in \overline{\rho }$
for every element $a\in I$. Hence, $(ae_{2},a)\in \theta $; similarly it can
be shown that $(e_{2}a,a)\in \theta $, and therefore, $\overline{e}_{2}$ is
the identity in the semiring $I/\theta $.

Suppose $K$ is a nonzero ideal in $I/\theta $, and let $\overline{K}:=\{x\in
I\,|\,\overline{x}\in K\}$. For each $a\in eRe$ and $x\in \overline{K}$, we
have $\overline{x}\in K$, $e_{2}a,ax\in I$; hence, $\overline{ax}=\overline{%
e_{2}}$ $\overline{ax}=\overline{e_{2}ax}=\overline{e_{2}a}$ $\overline{x}%
\in K$, that is, $ax\in \overline{K}$. Similarly, we also have $xa\in 
\overline{K}$. So, $\overline{K}$ is an ideal of $eRe$. Then, we have $R%
\overline{K}R$ is an ideal of $R$ and $R\overline{K}R\subseteq A$. Consider $%
L:=\{\overline{y}\,|\,y\in R\overline{K}R\}\subseteq A/\rho $. Clearly, $L$
is an ideal of $A/\rho $. Since $K\neq 0$, there is an element $k\in I$ such
that $\overline{0}\neq \overline{k}\in K$. This shows that $k\in \overline{K}
$ and $(e_{1}ke_{1},0)\notin \rho $, and hence, $\overline{0}\neq \overline{%
e_{1}ke_{1}}\in L$. For $A/\rho $ is ideal-simple, $L=A/\rho $; hence, $%
\overline{e}_{1}\in L$, \textit{i.e.}, there exist elements $%
r_{1},...,r_{l},s_{1},...,s_{l}\in R$ and $x_{1},...,x_{l}\in \overline{K}$
such that $\overline{e_{1}}=\overline{r_{1}x_{1}s_{1}+...+r_{l}x_{l}s_{l}}$.
As $x_{i}\in \overline{K}$, we have $ex_{i}e=x_{i}$ for all $i=1,\ldots l$.
Then, as $\overline{K}$ is an ideal of $eRe$, for all $i=1,\ldots l$, we
have $er_{i}x_{i}s_{i}e=er_{i}ex_{i}es_{i}e\in \overline{K}$ and $%
er_{1}x_{1}s_{1}e+...+er_{l}x_{l}s_{l}e\in \overline{K}$. Whence, $\overline{%
e}_{2}=\overline{ee_{1}e}=\overline{er_{1}x_{1}s_{1}e+...+er_{l}x_{l}s_{l}e}%
\in K$, and hence, $K=I/\theta $. Thus, $I/\theta $ is an ideal-simple
semiring, and therefore, $\mathcal{R}_{BM}(eRe)=0.$

(ii). As $\mathcal{R}_{BM}(R)$ is a subtractive ideal of $R$, the ideal $e%
\mathcal{R}_{BM}(R)e$ of $eRe$ is also subtractive. Similarly to the proof
of Proposition 5.15 (ii), it can be shown that $eRe/e\mathcal{R}%
_{BM}(R)e\cong \overline{e}$$\overline{R}$$\overline{e}$, where $\overline{e}
$ is the image of $e$ in $\overline{R}=R/\mathcal{R}_{BM}(R)$.

From $\mathcal{R}_{BM}(R/\mathcal{R}_{BM}(R))=0$ and (i), it follows that

\noindent $\mathcal{R}_{BM}(eRe/e\mathcal{R}_{BM}(R)e)=0$ and by %
\cite[Theorem 4.9]{m:otrtfs} we have the inclusion $\mathcal{R}%
_{BM}(eRe)=\cap \{K\in \mathcal{SI}(eRe)\,|\,\mathcal{R}_{BM}(eRe/K)=0\}%
\subseteq e\mathcal{R}_{BM}(R)e$. And we need only to show that the opposite
inclusion $e\mathcal{R}_{BM}(R)e\subseteq \mathcal{R}_{BM}(eRe)$ also takes
place.

So, suppose $e\mathcal{R}_{BM}(R)e\nsubseteq \mathcal{R}_{BM}(eRe)$. Then
there exists a congruence $\rho $ on $e\mathcal{R}_{BM}(R)e$ such that $e%
\mathcal{R}_{BM}(R)e/\rho $ is an ideal-simple semiring with identity $%
\overline{e_{1}}\,\ $with $e_{1}\in e\mathcal{R}_{BM}(R)e$. Similarly to the
proof in (i) and applying Lemma 3.14 and \cite[Propostion 5.3]{knz:ososacs},
we have that 
\begin{equation*}
\overline{\rho }=\{(a,b)\in (eRe)^{2}\,|\,(e_{1}ae_{1},e_{1}be_{1})\in \rho
\}
\end{equation*}%
is a congruence on $eRe$ and the relation $\theta $ on $R$, defined for all $%
a,b\in R$ by 
\begin{equation*}
(a,b)\in \theta \Longleftrightarrow \forall r,s\in R:(erase,erbse)\in 
\overline{\rho }\text{,}
\end{equation*}%
is a congruence on $R$. Now we shall show that for the congruence $\delta
:=\theta \cap \mathcal{R}_{BM}(R)^{2}$ on $\mathcal{R}_{BM}(R)$, the
hemiring $\mathcal{R}_{BM}(R)/\delta $ is, in fact, an ideal-simple
semiring. For $e$ is a full idempotent, $1=x_{1}ey_{1}+...+x_{n}ey_{n}$ for
some $x_{i},y_{i}\in R$, $i=1,\ldots ,n$, and $ee_{1}=e_{1}=e_{1}e\in e%
\mathcal{R}_{BM}(R)e\subseteq \mathcal{R}_{BM}(R)$. For any $x\in \mathcal{R}%
_{BM}(R)$ and $r,s\in R$, we have $%
e_{1}erx_{i}e_{1}y_{i}xsee_{1}=ee_{1}erx_{i}e.e_{1}.ey_{i}xsee_{1}e$ and $%
ee_{1}erx_{i}e,ey_{i}xsee_{1}e\in e\mathcal{R}_{BM}(R)e$. Whence, $%
(e_{1}erx_{i}e_{1}y_{i}xsee_{1},e_{1}erx_{i}ey_{i}xsee_{1})$

\noindent $\in \rho $ for all $i=1,...,n$, and $(e_{1}er(\Sigma
_{i=1}^{n}x_{i}e_{1}y_{i})xsee_{1},e_{1}erxsee_{1})=$

\noindent $\Sigma
_{i=1}^{n}(e_{1}erx_{i}e_{1}y_{i}xsee_{1},e_{1}erx_{i}ey_{i}xsee_{1})$ $\in
\rho $, and hence, $((\Sigma _{i=1}^{n}x_{i}e_{1}y_{i})x,x)\in \theta $.
Furthermore, for $e_{2}:=\Sigma _{i=1}^{n}x_{i}e_{1}y_{i}\in \mathcal{R}%
_{BM}(R)$ and any $x\in \mathcal{R}_{BM}(R)$, we have $(e_{2}x,x)\in 
\mathcal{R}_{BM}(R)^{2}$, and hence, $(e_{2}x,x)\in \delta $; similarly, $%
(xe_{2},x)\in \delta $ too. Thus, $\mathcal{R}_{BM}(R)/\delta $ is a
semiring with identity $\overline{e_{2}}$.

Let $I$ be a nonzero ideal of $\mathcal{R}_{BM}(R)/\delta $, and $\overline{I%
}:=\{x\in \mathcal{R}_{BM}(R)|\overline{x}\in I\}$. For each $r\in R$ and $%
x\in \overline{I}$, we have $e_{2}r,rx\in \mathcal{R}_{BM}(R)$, and hence, $%
\overline{rx}=\overline{e_{2}}\cdot \overline{rx}=\overline{e_{2}r}\cdot 
\overline{x}\in I$ and $rx\in \overline{I}$. Similarly, one also gets $xr\in 
\overline{I}$. Therefore, $\overline{I}$ is an ideal of $R$. It is clear
that $K:=\{\overline{x}$ $|$ $x\in e\overline{I}e\}\subseteq e\mathcal{R}%
_{BM}(R)e/\rho $ is an ideal of $e\mathcal{R}_{BM}(R)e/\rho $. For $I\neq 0$%
, there exists an element $a\in \mathcal{R}_{BM}(R)$ such that $0\neq 
\overline{a}\in I$. It shows that $a\in \overline{I}$ and $(a,0)\notin
\delta $, that is, there exist $r,s\in R$ such that $(erase,0)\notin 
\overline{\rho }$, \textit{i.e.}, $(e_{1}erasee_{1},0)\notin \rho $, and $%
K\neq 0$ because $e_{1}erasee_{1}\in K$. For $e\mathcal{R}_{BM}(R)e/\rho $
is ideal-simple, $K=e\mathcal{R}_{BM}(R)e/\rho $. Then, $\overline{e_{1}}\in
K$, that is, $e_{1}\in e\overline{I}e$. As $\overline{I}$ is an ideal of $R$%
, we have $e\overline{I}e\subseteq \overline{I}$, and hence, $e_{1}\in 
\overline{I}$ and $e_{2}=\Sigma _{i=1}^{n}x_{i}e_{1}y_{i}\in \overline{I}$.
Whence, $\overline{e_{2}}\in I$, and hence, $I=\mathcal{R}_{BM}(R)/\delta $.
Therefore, $\mathcal{R}_{BM}(R)/\delta $ is an ideal-simple semiring;
however, $\mathcal{R}_{BM}(R)$ has no nonzero semiring images. Thus, $e%
\mathcal{R}_{BM}(R)e\subseteq \mathcal{R}_{BM}(eRe)$.

(iii). Similarly to the proof of Proposition 5.1, we readily have that $P=Re$
is a progenerator of $_{R}\mathcal{M}$ and $P^{\ast
}=Hom_{R}(_{R}P,_{R}R)\cong eR$, $S=End_{R}(_{R}P)\cong eRe$, as well as
there exist the isomorphisms $\alpha :P\otimes _{eRe}P^{\ast }=Re\otimes
_{eRe}eR\longrightarrow R$ and $\beta :P^{\ast }\otimes _{R}P=eR\otimes
_{R}Re\longrightarrow eRe$ such that $\alpha (re\otimes er^{\prime
})=rer^{\prime }$ and $\beta (er\otimes r^{\prime }e)=err^{\prime }e$. Also,
similarly to as it was done in the proof of \cite[Corollary 4.4]{kn:meahcos}%
, we get the inverse category equivalences: $F:$ $_{R}\mathcal{M}%
\rightleftarrows $ $_{eRe}\mathcal{M}:G$ given by $F(A)=P^{\ast }\otimes
_{R}A=eR\otimes _{R}A\cong eA$ and $G(B)=P\otimes _{eRe}B=Re\otimes _{eRe}B$%
. Therefore, by Proposition 5.6, the functors $F$ and $G$ establish the
equivalences between the categories of simple semimodules of the categories $%
_{R}\mathcal{M}$ and $_{eRe}\mathcal{M}$, respectively.

Let $\mathcal{S}_{R}$ stay for the set of all simple left $R$-semimodules.
Then, $J_{s}(R)=\cap _{M\in \mathcal{S}_{R}}(0:M)_{R}$ and $J_{s}(eRe)=\cap
_{M\in \mathcal{S}_{R}}(0:F(M))_{eRe}$ $=\cap _{M\in \mathcal{S}%
_{R}}(0:eM)_{eRe}$. And we shall prove that $eJ_{s}(R)e=\cap _{M\in \mathcal{%
S}_{R}}(0:eM)_{eRe}$. Indeed, if $ere\in \cap _{M\in \mathcal{S}%
_{R}}(0:eM)_{eRe}$, then $ereM=ere.eM=0$ for all $M\in \mathcal{S}_{R}$, and
therefore, $ere\in J_{s}(R)$ and $ere\in eJ_{s}(R)e$, and hence, $\cap
_{M\in \mathcal{S}_{R}}(0:eM)_{eRe}\subseteq eJ_{s}(R)e$. Conversely, if $%
r\in J_{s}(R)$, then $rM=0$ for all $M\in \mathcal{S}_{R}$; so, $erM=0$ for
all $M\in \mathcal{S}_{R}$. As $eM\subseteq M$, we have $ere.eM=ereM$ for
all $M\in \mathcal{S}_{R}$. Whence, $ere\in \cap _{M\in \mathcal{S}%
_{R}}(0:eM)_{eRe}$, and hence, $eJ_{s}(R)e\subseteq \cap _{M\in \mathcal{S}%
_{R}}(0:eM)_{eRe}$.

The rest can be proved in the similar to the proof of Proposition 5.15 (ii)
way. \textit{\ \ \ \ \ \ }$_{\square }\medskip $

Our next result shows that the $J$-semisimplicity, $J_{s}$-semisimplicity,
and $\mathcal{R}_{BM}$-semisimplicity of semirings are Morita invariant
properties, namely:\medskip

\noindent \textbf{Theorem 5.17.} \textit{For Morita equivalent semirings} $R$
\textit{and} $S$\textit{, the followings statements are true:}

\textit{(i)} $J(R)=0$ \textit{iff }$J(S)=0$\textit{, and the semirings} $%
R/J(R)$ \textit{and} $S/J(S)$ \textit{are Morita equivalent;}

\textit{(ii) }$\mathcal{R}_{BM}(R)=0$ \textit{iff }$\mathcal{R}_{BM}(S)=0$%
\textit{,} \textit{and} \textit{the semirings} $R/\mathcal{R}_{BM}(R)$ 
\textit{and} $S/\mathcal{R}_{BM}(S)$ \textit{are Morita equivalent;}

\textit{(iii) } $J_{s}(R)=0$ \textit{iff } $J_{s}(S)=0$\textit{,} \textit{and%
} \textit{the semirings} $R/J_{s}(R)$ \textit{and} $S/J_{s}(S)$ \textit{are
Morita equivalent}\medskip .

\noindent \textbf{Proof.} (i). The assertion that $J(R)=0$ iff $J(S)=0$
follows right away from Proposition 5.1 and Theorem 5.8 or Corollary 5.11.

Next, by Proposition 5.1, $S\cong eM_{n}(R)e$ for some idempotent $e$ $\in
M_{n}(R)$ such that $M_{n}(R)eM_{n}(R)=M_{n}(R)$. Then, it is clear that $%
M_{n}(R/(J(R))\cong M_{n}(R)/M_{n}(J(R))$. From these observations,
Proposition 5.15 and Theorem 5.8, we have that $S/J(S)\cong \overline{e}%
M_{n}(R/J(R))\overline{e}$ and $M_{n}(R/J(R))\overline{e}M_{n}(R/J(R))$

\noindent $=M_{n}(R/J(R))$, and, using Proposition 5.1 again, we obtain that 
$R/J(R)$ and $S/J(S)$ are Morita equivalent.

Just using Proposition 5.16, (ii) and (iii) are proved in the same fashion
as (i). \textit{\ \ \ \ \ \ }$_{\square }\medskip $

Two hemirings $R$ and $S$, we say, are \emph{D}-\emph{Morita equivalent} if
their Dorroh extensions $R^{1}$ and $S^{1}$, respectively, are Morita
equivalent semirings. Since, by \cite[Corollary 4.4]{kn:meahcos}, the Morita
equivalence relation on the category of semirings is an equivalence
relation, the D-Morita equivalence relation for hemirings is also an
equivalence relation on the category of hemirings.\medskip

Using Lemmas 5.9, 5.10 and Theorem 5.17, we conclude the paper with the
following\medskip\ result:

\noindent \textbf{Corollary 5.18. }\textit{For D-Morita equivalent hemirings}
$R$ \textit{and} $S$\textit{, the followings statements are true:}

\textit{(i)} $J(R)=0$ \textit{iff }$J(S)=0$\textit{;}

\textit{(ii) }$\mathcal{R}_{BM}(R)=0$ \textit{iff }$\mathcal{R}_{BM}(S)=0$%
\textit{;}

\textit{(iii) } $J_{s}(R)=0$ \textit{iff } $J_{s}(S)=0$\medskip \textit{.}

\end{document}